\documentclass[journal]{IEEEtran}%
\usepackage{graphicx}
\usepackage{amsmath}
\usepackage{graphics}
\usepackage{amsfonts}
\usepackage{amssymb}%
\usepackage{listings}
\usepackage{amsthm}
\usepackage{epsfig}
\usepackage{psfrag}
\usepackage{booktabs}
\usepackage{hyperref}
\providecommand{\U}[1]{\protect\rule{.1in}{.1in}}

\newtheorem{thm}{Theorem}
\newtheorem{remark}{Remark}
\newtheorem{lemma}{Lemma}
\newtheorem{definition}{Definition}
\newtheorem{prop}{Proposition}
\newtheorem{coro}{Corollary}
\newtheorem{assum}{Assumption}
\newtheorem{property}{Property}
\usepackage{subfig}

\ifCLASSINFOpdf
\else
\fi
\begin{document}

\title{Sparse Regularization: Convergence Of Iterative Jumping Thresholding Algorithm}
\author{Jinshan Zeng, ~Shaobo Lin$^{\ast}$, ~and Zongben Xu
\thanks{ J.S. Zeng is with the School of Computer and Information Engineering, Jiangxi Normal University, Nanchang 330022,
and Beijing Center for
Mathematics and Information Interdisciplinary Sciences (BCMIIS),
Beijing, 100048, China. S.B. Lin is with the  College of Mathematics
and Information Science, Wenzhou University, Wenzhou 325035,
 Z.B. Xu is with the School of Mathematics and Statistics, Xi'an
Jiaotong University, Xi'an 710049, P R China. (email:
jsh.zeng@gmail.com, sblin1983@gmail.com, zbxu@mail.xjtu.edu.cn).
$*$ Corresponding author: Shaobo Lin (sblin1983@gmail.com). }}
\maketitle

\begin{abstract}
In recent studies on sparse modeling, non-convex penalties have received considerable attentions due to their superiorities on sparsity-inducing over the convex counterparts.
Compared with the convex optimization approaches, however,
the non-convex approaches have more challenging convergence analysis.
In this paper, we study the convergence of a non-convex iterative thresholding algorithm for solving sparse recovery problems with a certain class of non-convex penalties, whose corresponding thresholding functions are discontinuous with jump discontinuities.
Therefore, we call the algorithm the iterative jumping thresholding (IJT) algorithm.
The finite support and sign convergence of IJT algorithm is firstly verified via taking advantage of such jump discontinuity.
Together with the assumption of the introduced restricted Kurdyka-{\L}ojasiewicz (rKL) property,
then the strong convergence of IJT algorithm can be proved.
Furthermore, we can show that IJT algorithm converges to a local minimizer at an asymptotically linear rate under some additional conditions.
Moreover, we derive a posteriori computable error estimate, which can be used to design practical terminal rules for the algorithm.
It should be pointed out that the $l_q$ quasi-norm ($0<q<1$) is an important subclass of the class of non-convex penalties studied in this paper.
In particular, when applied to the $l_q$ regularization, IJT algorithm can converge to a local minimizer with an asymptotically linear rate under certain concentration conditions.
We provide also a set of simulations to support
the correctness of theoretical assertions and compare the time efficiency of
IJT algorithm for the $l_{q}$ regularization ($q=1/2, 2/3$)
with other known typical algorithms like the iterative reweighted least squares (IRLS) algorithm
and the iterative reweighted $l_{1}$ minimization (IRL1) algorithm.
\newline

\end{abstract}




\markboth{ ~}{Shell \MakeLowercase{\textit{et al.}}: Bare Demo of IEEEtran.cls for Journals}






\IEEEpeerreviewmaketitle

\section{Introduction}



The sparse vector recovery problems emerging in many areas of
scientific research and engineering practice have attracted
considerable attention in recent years
(\cite{Donoho06}-\cite{Duarte2011StructuredCS}). Typical
applications include regression {\cite{TibshiraniVS}}, visual coding
{\cite{OlshausenVC}}, signal processing
{\cite{Combettes2005ProxInSP}}, compressed sensing \cite{Donoho06},
{\cite{Candes06}}, machine learning {\cite{1-SVMZhu2003}}, and
microwave imaging {\cite{SARZeng2012}}.
These problems can be modeled as the following $l_0$-norm regularized optimization problem
\begin{equation}
\min_{x\in \mathbf{R}^N} \left\{ F(x) + \lambda \|x\|_0 \right\},
\label{L0Reg}
\end{equation}
where $F: \mathbf{R}^N \rightarrow [0,\infty)$ is a proper lower-semicontinuous function,
$\|x\|_0$, commonly called the $l_0$-norm, denotes the number of nonzero components of $x$ and $\lambda>0$ is a regularization parameter.
The $l_0$ regularized least squares problem is a special case of  (\ref{L0Reg}) where $F(x) = \frac{1}{2}\|Ax-y\|_2^2$.
Blumensath and Davies {\cite{Blumensath08}} proposed the iterative \textit{hard} thresholding algorithm to solve this problem,
and showed that the algorithm converges to a local minimizer.
Recently, Lu and Zhang {\cite{LU2013PDM}} proposed a penalty decomposition method for solving a more general class of $l_0$ regularized problems.
In addition, Lu {\cite{LU2014IHT}} proposed an iterative \textit{hard} thresholding method and its variant for solving $l_0$ regularization over a conic constraint, and established its convergence as well as the iteration complexity.


Besides the $l_0$ regularized optimization problem, a more general class of problems are considered a lot in both practice and theory, that is,
\begin{equation}
\min_{x\in \mathbf{R}^N} \{F(x) + \lambda \Phi(x)\},
\label{NonSCS}
\end{equation}
where $\Phi(x)$ is a certain separable, continuous penalty with $\Phi(x) = \sum_{i=1}^N \phi(|x_i|)$, and $x=(x_1,\cdots,x_N)^T$.
One of the most important cases is the $l_1$-norm with $\Phi(x) = \|x\|_1=\sum_{i=1}^N |x_i|$.
The $l_1$-norm is convex and thus, the corresponding $l_1$-norm regularized optimization problem can be efficiently solved.
Because of this, the $l_1$-norm becomes popular and has been accepted as a very useful tool for the modeling of the sparsity problems.
Nevertheless, the $l_1$-norm may not induce adequate sparsity when applied to certain applications {\cite{Candes2008RL1}}, {\cite{Chartrand2007}}, {\cite{Chartrand2008}}, {\cite{L1/2TNN}}.
Alternatively, many non-convex penalties were proposed as relaxations of the $l_0$-norm.
Some typical non-convex examples are the $l_q$-norm ($0<q<1$) {\cite{Chartrand2007}}, {\cite{Chartrand2008}}, {\cite{L1/2TNN}}, Smoothly Clipped Absolute Deviation (SCAD) {\cite{FanSCAD}}, Minimax Concave Penalty (MCP) {\cite{ZhangMCP2010}} and Log-Sum Penalty (LSP) {\cite{Candes2008RL1}}.
Compared with the $l_1$-norm, the non-convex penalties can usually induce better sparsity while the corresponding non-convex regularized optimization problems are generally more difficult to solve.

There are mainly four classes of algorithms to solve the non-convex
regularized optimization problem (\ref{NonSCS}). The first one is
the half-quadratic (HQ) algorithm {\cite{GRHQ1992}},
{\cite{GYHQ1995}}.
HQ algorithms
can be efficient when both subproblems are easy to solve
(particularly, when both subproblems have closed-form solutions).
The second class is the iterative reweighted algorithm including iterative reweighted least squares (IRLS) minimization
{\cite{FOCUSS}}, {\cite{IRAChartrand_Yin2008}},
{\cite{DaubechiesIRLS}}  and iterative reweighted $l_1$-minimization
(IRL1) {\cite{Candes2008RL1}} algorithms.
Recently, Lu {\cite{LU2014IRM}} extended some existing
iterative reweighted methods and then proposed new variants for the
general $l_q$ ($0<q<1$) regularized unconstrained minimization
problems. Nevertheless, the iterative reweighted algorithms can be
only efficient when
the corresponding non-convex penalty can be well approximated via the quadratic function or the weighted $l_1$-norm function.
The third class is the difference of convex functions algorithm (DC programming) {\cite{DCProgGasso2009}}, which is also called
Multi-Stage (MS) convex relaxation {\cite{ZhangMS2010}}.
The DC programming considers a proper decomposition of the objective function.
Hence, it can be only applied to those non-convex penalties that can be decomposed as a difference of convex functions.
The last class is the iterative thresholding algorithm, which fits the framework of the forward-backward splitting algorithm {\cite{Attouch2013}} and the generalized gradient projection algorithm {\cite{BrediesNonconvex}} when applied to a separable non-convex penalty.
Intuitively, the iterative thresholding algorithm can be viewed as a procedure of
Landweber iteration projected by a certain thresholding operator.
Thus, the thresholding operator plays a key role in the iterative thresholding algorithm.
For some special non-convex penalties such as SCAD, MCP, LSP and $l_q$-norms with $q=1/2, 2/3$,
the associated thresholding operators can be expressed analytically {\cite{L1/2TNN}}, {\cite{L2/3Cao2013}}, {\cite{YeNonconvex}}.
Compared to the other types of non-convex algorithms such as the HQ, IRLS, IRL1 and DC programming algorithms,
the iterative thresholding algorithm is easy to implement and
has almost the least computational complexity for large scale problems {\cite{SARZeng2012}}, {\cite{Qian20011}}.
Consequently, the iterative thresholding algorithm becomes popular.

One of the significant differences between the convex and non-convex algorithms is that the convergence analysis of a non-convex algorithm is
in general tricky.
Although the effectiveness of the iterative thresholding algorithms for the non-convex regularized optimization problems has been verified in many
applications, except for the iterative \textit{hard} {\cite{LU2014IHT}} and \textit{half} {\cite{ZengHalfConv2013}} thresholding algorithms,
the convergence of most of these algorithms has not been thoroughly
investigated. More specifically, there are still three mainly open
questions.
\begin{enumerate}
  \item When does the algorithm converge? Under what conditions, the iterative thresholding algorithm converges strongly
  in the sense that the whole sequence generated, regardless of the initial point, is convergent.
  \item Where does the algorithm converge? Does the algorithm converge to a global minimizer or more practically, a local minimizer
due to the non-convexity of the optimization problem?
  \item What is the convergence rate of the algorithm?
\end{enumerate}

\subsection{Main Contribution}

In this paper, we give the convergence analysis for the iterative jumping thresholding  algorithm (called IJT algorithm henceforth) for solving a certain class of non-convex regularized optimization problems.
One of the most significant features of such non-convex problems is that the corresponding thresholding functions are discontinuous with jump discontinuities (see Fig. {\ref{Fig_ThresFun}}).
Moreover, the corresponding thresholding functions are not nonexpansive in general.
Among these non-convex penalties, the well-known $l_q$-norm with $0<q<1$
is one of the most typical cases. The main contribution can be summarized as follows.

\begin{enumerate}
\item[(a)]
We prove that the supports and signs of any sequence generated by IJT algorithm can converge within finite iterations.
Such property brings a possible way to construct a new sequence in a special subspace such that the new sequence has the same convergence behavior of the original sequence generated by IJT algorithm.

\item[(b)]
Under a further assumption that the objective function satisfies the so-called restricted Kurdyka-{\L}ojasiewicz (rKL) property (see Definition {\ref{Def_rKLProp}}) at some limit point, the strong convergence of IJT algorithm can be assuredly guaranteed (see Theorem {\ref{Thm_StrongConv}}).
The introduced rKL property is generally weaker than the well-known Kurdyka-{\L}ojasiewicz property that is widely used to study the convergence of nonconvex algorithms.

\item[(c)]
Under certain second-order conditions, we demonstrate that IJT algorithm converges to a local minimizer at an asymptotically linear rate (see Theorems {\ref{Thm_LocalConv}}-{\ref{Thm_ConvRate2}}).
Such asymptotically linear convergence speed means that when the iterative vector is sufficiently close to the convergent point, the rate of convergence
of IJT algorithm is linear. This implies that given a good initial guess, IJT algorithm can converge very fast.

\item[(d)]
As a typical case, we apply the developed convergence results to the $l_q$ regularization ($0<q<1$).
When applied to the $l_q$ regularization, IJT algorithm can converge to a local minimizer at an asymptotically linear rate
as long as the matrix satisfies a certain concentration property (see Theorem {\ref{Thm_ConvRateLq_mu}}).

\item[(e)]
We also provide simulations to
support the correctness of theoretical assertions and compare the
convergence speed of IJT algorithm for the $l_{q}$
regularization problems ($q=1/2, 2/3$) with other known typical algorithms
like the iterative reweighted least squares (IRLS) algorithm and the
iterative reweighted $l_{1}$ minimization (IRL1) algorithm.
\end{enumerate}

\subsection{Notations and Organization}

We denote $\mathbf{R}$ and $\mathbf{N}$ as the real number and natural number sets, respectively.
For any vector $x\in \mathbf{R}^N$, $x_i$ is its $i$-th component,
and for a given index set $I\subset I_N \triangleq \{1,2,\cdots,N\}$, $x_I$ represents its subvector containing all the components restricted to $I$.
$I^c$ represents the complementary set of $I$, i.e., $I^c = I_N \setminus I.$
$\|x\|_2$ represents the Euclidean norm of a vector $x$.
$Supp(x)$ is the support set of $x$, i.e., $Supp(x) = \{i: |x_i|>0, i=1,\cdots, N\}$.
For any matrix $A\in \mathbf{R}^{N\times N}$,
$\sigma_i(A)$ and $\sigma_{\min}(A)$ ($\lambda_i(A)$ and $\lambda_{\min}(A)$) denote as the $i$-th and minimal singular values (eigenvalues) of $A$, respectively.
Similar to the vector case, for a given index set $I$, $A_I$ represents the submatrix of $A$ containing all the columns restricted to $I$.
For any $z\in \mathbf{R}$, $sign(z)$ denotes its sign function, i.e.,
\[
sign(z)=\left\{
\begin{array}
[c]{ll}%
1, & \mbox{for} \ z>0\\
0, & \mbox{for} \ z=0\\
-1, & \mbox{for} \ z<0\\
\end{array}
\right..
\]

The remainder of this paper is organized as follows.
In section II, we give the problem settings and then introduce IJT algorithm with some basic properties.
In section III, we give the convergence analysis of IJT algorithm.
In section IV, we apply the established theoretical analysis to the $l_q$ ($0<q<1$) regularization.
In section V, we discuss some related work.
In section VI, we conduct the simulations to substantiate the theoretical results.
We conclude this paper in section VII.

\section{Iterative Jumping Thresholding Algorithm}

In this section, we first present the basic settings of the considered non-convex regularized optimization problems,
then introduce IJT algorithm for these problems.
In the end of this section, we briefly review some basic properties of IJT algorithm obtained in {\cite{BrediesNonconvex}}.

\subsection{Problem Settings}

We consider the following composite optimization problem
\begin{equation}
\min_{x\in \mathbf{R}^N} \{T_{\lambda}(x) = F(x) + \lambda \Phi(x)\},
\label{NonOpt}
\end{equation}
where $\Phi(x)$ is assumed to be separable with $\Phi(x) = \sum_{i=1}^N \phi(|x_i|)$.
Moreover, we make several assumptions on the problem (\ref{NonOpt}).

\begin{assum} \label{AonF}
$F: \mathbf{R}^N \rightarrow [0,\infty)$ is weakly lower-semicontinuous and differentiable
with Lipschitz continuous gradient, i.e., it holds that
$$
\|\nabla F(u)- \nabla F(v)\|_2 \leq L\|u-v\|_2,\ \  \forall u, v \in \mathbf{R}^N,
$$
where $L>0$ is the Lipschitz constant.
\end{assum}
%

It should be noted that Assumption {\ref{AonF}} is a general assumption for
$F$. Many formulations in machine learning satisfy Assumption {\ref{AonF}}. For
example, the following least squares and logistic loss functions are
two commonly used functions which satisfy Assumption {\ref{AonF}}:
$$F(x) = \frac{1}{2M}\|Ux-y\|_2^2 \ \ or\ \  \frac{1}{M}\sum_{i=1}^M \log(1+ \exp(-y_iu_i^Tx)),$$
where $u_i \in \mathbf{R}^N$ for $i=1,2,\cdots, M$, $U=[u_1, \cdots, u_M]^T \in \mathbf{R}^{M\times N}$ is a data matrix and $y=(y_1,\cdots,y_M)^T \in \mathbf{R}^M$ is a target vector.
Moreover, in both signal and image processing, $F$ is commonly taken as the least squares of the observation model, that is,
$$F(x) = \|Ax-y\|_2^2,$$
where $y\in \mathbf{R}^M$ is an observation vector and $A\in \mathbf{R}^{M\times N}$ is an observation matrix.
It can be easily verified that such $F$ also satisfies Assumption {\ref{AonF}}.

In the following, we give some basic assumptions on $\phi$, most of which were considered in {\cite{BrediesNonconvex}}.

\begin{assum}\label{AonPhi}
$\phi : [0,\infty) \rightarrow [0,\infty)$ is continuous and satisfies the following assumptions:

\begin{enumerate}
\item[(a)]
$\phi$ is non-decreasing with $\phi(0) = 0$ and $\phi(z)\rightarrow \infty$ when $z \rightarrow \infty$.

\item[(b)]
For each $b>0$, there exists an $a>0$ such that $\phi(z) \geq az^2$ for $z \in [0,b]$.

\item[(c)]
$\phi$ is differentiable on $(0,\infty)$ and the derivative $\phi'$ is strictly convex with $\phi'(z)\rightarrow \infty$ for $z\rightarrow 0$ and $\phi'(z)/z\rightarrow 0$ for $z \rightarrow \infty$.

\item[(d)]
$\phi$ has a continuous second derivative $\phi''$ on $(0,\infty)$.

\end{enumerate}
\end{assum}

\begin{figure}[!t]
\begin{minipage}[b]{.49\linewidth}
\centering
\includegraphics*[scale=0.33]{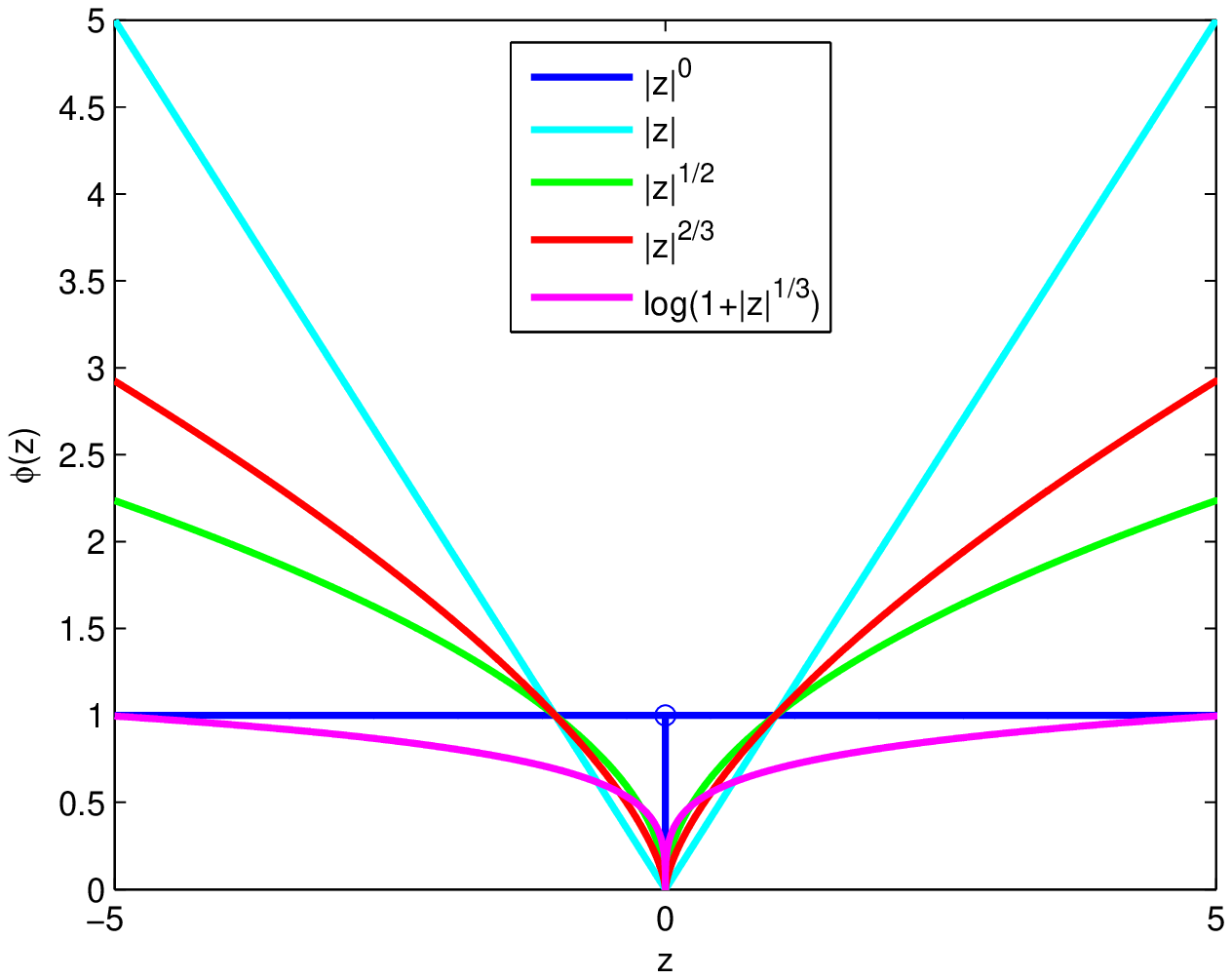}
\centerline{{\small (a) Typical Penalty Functions}}
\end{minipage}
\hfill
\begin{minipage}[b]{.49\linewidth}
\centering
\includegraphics*[scale=0.33]{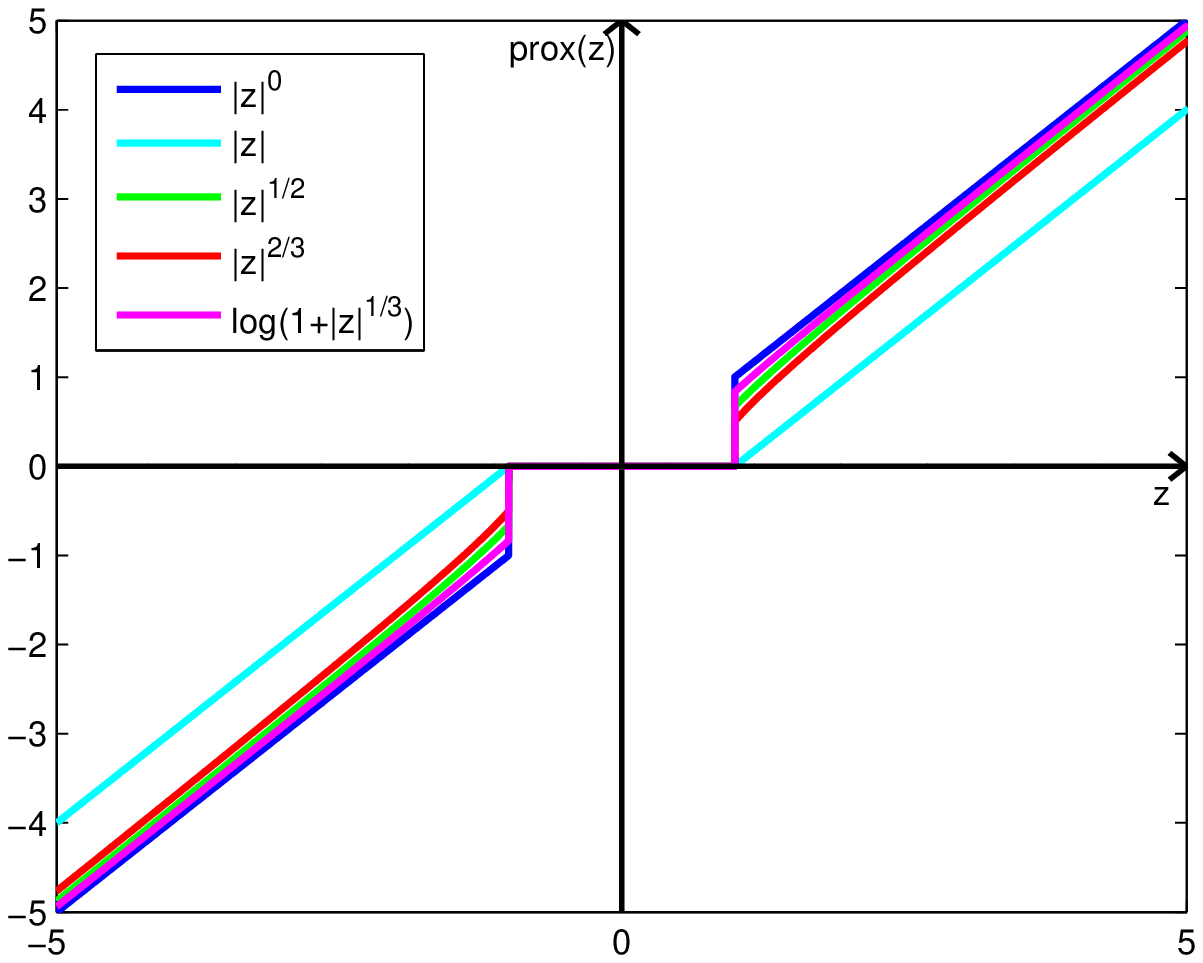}
\centerline{{\small (b) Thresholding Functions}}
\end{minipage}
\hfill
\caption{ Typical penalty functions $\phi$ satisfying Assumption {\ref{AonPhi}} and the corresponding thresholding functions.
More specifically, we plot the figures of the penalty functions $\phi(|z|) = |z|^{1/2}, |z|^{2/3}, \log(1+|z|^{1/3})$, and their corresponding thresholding functions.
For comparison, we also plot the figures of two well-known cases, i.e., $l_0$-norm with $\phi(|z|) = \textsl{1}_{|z|>0}$ as the indicator function of $|z|>0$, $l_1$-norm with $\phi(|z|)=|z|$, and their corresponding thresholding functions.
(a) Typical penalty functions.
(b) Thresholding functions.
}
\label{Fig_ThresFun}
\end{figure}

In Assumption {\ref{AonPhi}}, (a) and (b) are taken from Assumption 3.1 in {\cite{BrediesNonconvex}}, while (c) and (d) are adapted from Assumption 3.2 in {\cite{BrediesNonconvex}}.
It can be observed that Assumption {\ref{AonPhi}}(a) ensures the coercivity of $\phi$, and thus the existence of the minimizer of the optimization problem (\ref{NonOpt}).
Assumption {\ref{AonPhi}}(b) guarantees the weakly sequential lower semi-continuity of $\phi$ in $\l^2$,
and Assumption {\ref{AonPhi}}(c) induces the sparsity of the penalty $\Phi$.
In practice, there are many non-convex functions satisfying Assumption {\ref{AonPhi}}.
Two of the most typical subclasses are $\phi(z) = z^q$ and $\phi(z) = \log(1+z^q)$ with $q\in (0,1)$ as shown in Fig. {\ref{Fig_ThresFun}}.

\subsection{IJT Algorithm}

In order to describe IJT algorithm,
we need to generalize the proximity operator from the convex case to a non-convex penalty $\Phi$, that is,
\begin{equation}
Prox_{\mu,\lambda\Phi}(x) = \arg \min_{u\in \mathbf{R}^N} \left\{\frac{\|x-u\|_2^2}{2\mu} + \lambda \Phi(u)\right\},
\label{ProxOper}
\end{equation}
where $\mu>0$ is a parameter. Since $\Phi$ is separable,  computing
$Prox_{\mu,\lambda\Phi}$ is reduced to solve a one-dimensional
minimization problem, that is,
\begin{equation}
prox_{\mu,\lambda\phi}(z) = \arg \min_{v\in \mathbf{R}} \left\{\frac{|z-v|^2}{2\mu} + \lambda \phi(|v|)\right\}.
\label{SVProxOper}
\end{equation}
Therefore,
\begin{equation}
Prox_{\mu,\lambda\Phi}(x) = (prox_{\mu,\lambda\phi}(x_1), \cdots, prox_{\mu,\lambda\phi}(x_N))^T.
\end{equation}

As shown by (\ref{SVProxOper}), the proximity operator is defined through an optimization problem,
which is commonly hard for computing and analysis.
In order to present a simpler form of the proximity operator for analysis, we show a preparatory lemma in the following.

\begin{lemma} \label{Lemm_Pre}
{\bf (Lemma 3.10 in {\cite{BrediesNonconvex}})}
Assume that $\phi$ satisfies Assumption {\ref{AonPhi}}, then
\begin{enumerate}
\item[(a)]
for each $\mu>0$, the function $\rho_{\mu}: z \mapsto z + \lambda \mu \phi'(z)$ is well defined on $\mathbf{R}_{+}$ and,
moreover, it is strictly convex and attains a minimum at $z_{\mu}>0$;

\item[(b)]
the function $\psi: z \mapsto 2(\phi(z)-z\phi'(z))/z^2$ is strictly decreasing and one-to-one on $(0,\infty) \rightarrow (0,\infty)$;

\item[(c)]
for any $z>0$, it holds that $\phi''(z)<-\psi(z)<0$;

\item[(d)]
for any $z>0$, $\phi''(z)$ is negative and monotonically increasing.
\end{enumerate}
\end{lemma}

With Lemma {\ref{Lemm_Pre}}, $prox_{\mu,\lambda\phi}$ can be expressed as follows.

\begin{lemma}\label{Lemm_JumpThresFun}
{\bf (Lemma 3.12 in {\cite{BrediesNonconvex}})}
Assume that $\phi$ satisfies Assumption 2,
then $prox_{\mu, \lambda\phi}$ is well defined and can be specified as
\begin{equation}
prox_{\mu, \lambda\phi} (z)=\left\{
\begin{array}
[c]{ll}%
sign(z)\rho_{\mu}^{-1}(|z|), & \mbox{for} \ |z|\geq \tau_{\mu}\\
0, & \mbox{for} \ |z|\leq \tau_{\mu}
\end{array}
\right.,
\label{ProxMapExp}
\end{equation}
for any $z \in \mathbf{R}$ with
\begin{equation}
\tau_{\mu} = \rho_{\mu}(\eta_{\mu})
\label{ThreshValuex}
\end{equation}
and
\begin{equation}
\eta_{\mu} = \psi^{-1}((\lambda\mu)^{-1}).
\label{ThreshValuey}
\end{equation}
Moreover, the range of $prox_{\mu,\lambda\phi}$ is $\{0\}\cup [\eta_{\mu},\infty)$.
\end{lemma}

It can be observed that the proximity operator is discontinuous with a jump discontinuity,
which is one of the most significant features of such a class of non-convex penalties studied in this paper.
Moreover, it can be easily checked that the proximity operator is not nonexpansive in general.
Due to these, the convergence analysis of the corresponding non-convex algorithm gets challenging.
(Some specific proximity operators are shown in Fig. {\ref{Fig_ThresFun}}(b).)

With the definition of the proximity operator, IJT algorithm can be proposed to solve the non-convex regularized optimization problem (\ref{NonOpt}).
Formally, the iterative form of IJT algorithm can be expressed as follows
\begin{equation}
x^{n+1} \in Prox_{\mu, \lambda\Phi}(x^n-\mu \nabla F(x^n)),
\label{GGPA}
\end{equation}
where $\mu>0$ is a step size parameter.
For simplicity, we define
$$
G_{\mu,\lambda\Phi}(x) = Prox_{\mu, \lambda\Phi}(x-\mu \nabla F(x))
$$
for any $x\in \mathbf{R}^N$.
Henceforth, we call $prox_{\mu,\lambda\phi}$ the jumping thresholding function.

\begin{remark} \label{RThresq}
For some specific $l_q$-norm (say, $q=1/2, 2/3$), the proximity operator can be expressed analytically {\cite{L1/2TNN}}, {\cite{L2/3Cao2013}} (as shown in Fig. 1(b)).
\end{remark}

\begin{remark}\label{RHard}
Although the $l_0$-norm does not satisfy Assumption 2, the \textit{hard} thresholding function is also discontinuous with jump discontinuities.
Due to such discontinuity of the \textit{hard} thresholding function, we will discuss that the convergence of the \textit{hard} algorithm can be easily developed according to a similar analysis of IJT algorithm in Section III.
\end{remark}

\subsection{Some Basic Properties of IJT Algorithm}

In this subsection, we briefly review some basic properties of IJT algorithm, which serve as the basis of the further analysis in the next sections.
Some of these properties can be found in {\cite{BrediesNonconvex}}.

\begin{property}\label{Prop_SuffDecrease}
{\bf (Proposition 2.1 and Corollary 2.2 in {\cite{BrediesNonconvex}})} Let $\{x^n\}$ be a sequence generated by IJT algorithm
with a bounded initialization.
Assume that $0<\mu<\frac{1}{L}$, then it holds
\begin{enumerate}
\item[(a)]
$T_{\lambda}(x^{n+1}) \leq T_{\lambda}(x^n) - \frac{1}{2}(\frac{1}{\mu}-L) \|x^{n+1}-x^n\|_2^2$,
and there exists a positive constant $T_{\lambda}^*$ such that $T_{\lambda}(x^n) \rightarrow T_{\lambda}^*$ as $n \rightarrow \infty$;

\item[(b)]
$\|x^{n+1} - x^n\|_2 \rightarrow 0$ as $n \rightarrow \infty$.
\end{enumerate}
\end{property}

Property \ref{Prop_SuffDecrease}(a) is commonly called the sufficient decrease property, which is a basic property desired for a descent method.
With Property \ref{Prop_SuffDecrease}, the subsequential convergence of IJT algorithm can be easily claimed as the following property.

\begin{property}\label{Prop_SubseqConv}
{\bf (Proposition 2.3 in {\cite{BrediesNonconvex}}). }
Let $\{x^n\}$ be a sequence generated by IJT algorithm with a bounded initialization.
Suppose that $0<\mu<\frac{1}{L}$, then
\begin{enumerate}
\item[(a)]
each minimizer of $T_{\lambda}$ is a fixed point of $G_{\lambda\mu,\Phi}$;

\item[(b)]
there exists a convergent subsequence of $\{x^n\}$ and the limit point is a fixed point of $G_{\lambda\mu,\Phi}$.
\end{enumerate}
\end{property}

Besides Properties {\ref{Prop_SuffDecrease}} and \ref{Prop_SubseqConv}, we can derive the following property directly from the definition of the proximity operator.

\begin{property}\label{Prop_OptCond}
Let $x^*$ be a fixed point of $G_{\lambda\mu,\Phi}$ and $\{x^n\}$ be a sequence generated by IJT algorithm,
then it holds
\begin{enumerate}
\item[(a)]
$|x^*_i|\geq \tau_{\mu}/\mu$ and $[\nabla F(x^*)]_i + \lambda sign(x_i^*)\phi'(|x_i^*|)=0$ for any $i\in Supp(x^*)$, and $|[\nabla F(x^*)]_i|\leq \tau_{\mu}/\mu$ for any $i\in Supp(x^*)^c$;

\item[(b)]
$x_i^{n+1} + \lambda \mu sign(x_i^{n+1}) \phi'(|x_i^{n+1}|) = x_i^n - \mu [\nabla F(x^n)]_i$ for any $i\in Supp(x^{n+1})$
and $|x_i^n - \mu [\nabla F(x^n)]_i| \leq \tau_{\mu}$ for any $i\in Supp(x^{n+1})^c$,
$n \in \mathbf{N}$,
\end{enumerate}
where
$[\nabla F(x^*)]_i$ and $[\nabla F(x^{n+1})]_i$ represent the $i$-th component of $\nabla F(x^*)$ and $\nabla F(x^{n+1})$ respectively.
\end{property}

Actually, Property \ref{Prop_OptCond}(a) is a certain type of optimality conditions of the non-convex regularized optimization problem (\ref{NonOpt}).
Moreover, we call $x^*$ a \textit{stationary point} of (\ref{NonOpt}) if $x^*$ satisfies Property \ref{Prop_OptCond}(a),
and we denote $\Omega_{\mu}$ the stationary point for a given $\mu$.

\section{Convergence Analysis}

In the last section,  it can be only claimed that any sequence
$\{x^n\}$ generated by IJT algorithm subsequentially converges
to a stationary point. In this section, we will answer the open
questions concerning IJT algorithm presented in the
introduction, i.e., when, where and how fast does the algorithm
converge? More specifically, we first prove that IJT algorithm
converges to a stationary point under the so-called restricted Kurdyka-{\L}ojasiewicz (rKL) property (see Definition {\ref{Def_rKLProp}}), and
then show that the stationary point is also a local minimizer of the
optimization problem with some additional assumptions, and further
demonstrate that the convergence rate of IJT algorithm is
asymptotically  linear.

\subsection{Restricted Kurdyka-{\L}ojasiewicz Property}

Kurdyka-{\L}ojasiewicz (KL) property has been widely used to prove the convergence of the nonconvex algorithms (see, {\cite{Attouch2013}} for an instance).
Specifically, the KL property is the following.

\begin{definition}\label{Def_KLProp}
{\bf ({\cite{Attouch2013}})}
The function $f:\mathbf{R}^N \rightarrow \mathbf{R}\cup \{+\infty\}$ is said to have the Kurdyka-{\L}ojasiewicz property at $x^*\in$ dom $\partial f$ if there exist $\eta \in (0,+\infty]$, a neighborhood $U$ of $x^*$ and a continuous concave function $\varphi:[0,\eta)\rightarrow \mathbf{R}_{+}$ such that:
    \begin{enumerate}
    \item[(i)]
    $\varphi(0) = 0$;

    \item[(ii)]
    $\varphi$ is ${\cal{C}}^1$ on $(0,\eta)$;

    \item[(iii)]
    for all $s\in (0,\eta)$, $\varphi'(s)>0$;

    \item[(iv)]
    for all $x$ in $U\cap \{x: f(x^*)<f(x)<f(x^*) + \eta\}$,
    the Kurdyka-{\L}ojasiewicz inequality holds
    \begin{equation}
    \varphi'(f(x)-f(x^*)) dist(0,\partial f(x)) \geq 1.
    \label{KLIneq}
    \end{equation}
    \end{enumerate}
Proper lower semi-continuous functions which satisfy the Kurdyka-{\L}ojasiewicz inequality at each point of dom $\partial f$ are called KL functions.

\end{definition}

Roughly speaking, KL inequality means that the function considered is sharp up to a reparametrization at a neighborhood of some point.
From Definition {\ref{Def_KLProp}}, we can observe that KL inequality is actually certain type of first-order condition,
which implies that the gradient (subgradient or subdifferential)
of the transformed function via a concave function $\varphi$ is sharp and far away from zero.
Functions satisfying the KL inequality include real analytic functions, semialgebraic functions and locally strongly convex functions (more information can be referred to Sec. 2.2 in {\cite{Yin2013}} and references therein).

If further the objective function $T_{\lambda}$ in (\ref{NonOpt}) is a KL function
and the so-called relative error condition holds for the sequence $\{x^n\}$ generated by IJT algorithm,
then according to Theorem 5.1 in {\cite{Attouch2013}},
the strong convergence of IJT algorithm can naturally hold.
However, on one hand, the relative error condition may be violated for $\{x^n\}$.
Actually, as justified in the consequent Lemma {\ref{Lemm_NewSeq3Cond}},
such relative error condition only holds for the support sequence of $\{x^n\}$.
On the other hand,
as listed in Appendix A, we can construct a one-dimensional function that satisfies Assumptions {\ref{AonF}} and {\ref{AonPhi}}, but is not a KL function.
This motivates us to introduce the following so-called restricted Kurdyka-{\L}ojasiewicz (rKL) property to derive the convergence of IJT algorithm.
To describe the definition of rKL property conveniently, we define a projection mapping associated with an index set $I \subset \{1,2,\cdots,N\}$, that is,
$$
P_I: \mathbf{R}^N \rightarrow \mathbf{R}^K, P_Ix = x_I,  \forall x\in \mathbf{R}^N.
$$
We also denote $P_I^T$ as the transpose of $P_I$, i.e.,
$$
P_I^T: \mathbf{R}^{|I|} \rightarrow \mathbf{R}^N, (P_I^Tz)_I = z \ \text{and} \ (P_I^Tz)_{I^c} = 0, \forall z\in \mathbf{R}^{|I|},
$$
where $|I|$ is the cardinality of $I$ and $I^c = \{1,2,\cdots,N\} \setminus I$.

\begin{definition}\label{Def_rKLProp}
A function $f: \mathbf{R}^N \rightarrow \mathbf{R}\cup \{+\infty\}$ is said to have the $I$-restricted Kurdyka-{\L}ojasiewicz property at $x^*\in$ dom $\partial f$ with $I$ being a given subset of $\{1,2,\cdots, N\}$,
if the function $g: \mathbf{R}^{|I|} \rightarrow \mathbf{R}\cup \{+\infty\}, g(z)= f(P_I^T z)$ satisfies the KL inequality at $z^* = x_I^*.$
\end{definition}

Obviously, the introduced rKL property is weaker than the KL property.
If $I=\{1,2,\cdots,N\}$, then rKL property is exactly equivalent to the KL property.
From Definition {\ref{Def_rKLProp}}, rKL property only requires the subdifferential of the function with respect to a part of variables can get sharp after certain a concave transform, while KL property requires such well property for all the variables around some point.
It can be observed that rKL property is a natural extension of KL property.
Assume that $f_1: \mathbf{R}^{n_1} \rightarrow \mathbf{R}$ is a KL function, and $f_2: \mathbf{R}^{n_2} \rightarrow \mathbf{R}$ is an arbitrary function.
Let $f: \mathbf{R}^{n_1+n_2} \rightarrow \mathbf{R}, f(u)=f_1(u_{I_{n_1}}) + f_2(u_{I_{n_1}^c})$,
where $I_{n_1} = \{1, \cdots, n_1\}$ and $I_{n_1}^c = \{n_1+1, \cdots, n_1 + n_2\}$.
Then obviously, $f$ is a $I_{n_1}$-rKL function, but not a KL function.
In the following, we will give a sufficient condition of the rKL property.

\begin{lemma}
\label{Lemm_SuffCond_rKL}
Given an index set $I \subset \{1,2,\cdots,N\}$, consider the function $g(z)=f(P_I^Tz)$.
Assume that $z^*$ is a stationary point of $g$, and $g$ is twice continuously differentiable at a neighborhood of $z^*$, i.e., $B(z^*,\epsilon_0)$.
Moreover, if $\nabla^2 g(z^*)$ is nonsingular, then $f$ satisfies $I$-rKL property at the point $P_I^Tz^*$.
Actually, it holds
\[
|g(z) - g(z^*)| \leq C^* \|\nabla g(z)\|_2^2, \forall z\in B(z^*,\epsilon),
\]
for some $0<\epsilon<\epsilon_0$ and a positive constant $C^*>0$.
\end{lemma}

The proof of this lemma is shown in Appendix B.
From Lemma {\ref{Lemm_SuffCond_rKL}}, $g$ actually satisfies the KL inequality at $z^*$ with a desingularizing function of the form $\varphi(s)=c\sqrt{s},$ where $c>0$ is a constant.
Distinguished with the well-known KL inequality condition, the sufficient condition listed in the above lemma is some type of second-order condition, i.e., the Hessian of $g$ is nonsingular at some stationary point $z^*$.
The similar condition is also used to guarantee the convergence of the steepest descent method in {\cite{Ostrowski1967}} (Theorem 2, pp. 266).
Obviously, if a stationary point $z^*$ is a strictly local minimizer (or maximizer), or a strict saddle point of $g$,
then the nonsingularity of $\nabla^2 g(z^*)$ holds naturally.

\subsection{Convergence To A Stationary Point}

As analyzed in the section II, we have known that the sequence $\{x^n\}$ converges weakly.
Let $\cal{X}$ be the limit point set of $\{x^n\}$, $I^n = Supp(x^n)$.
In the following, we first show that both the support and sign of the sequence will converge within finite iterations,
and also any limit point $x^* \in \cal{X}$ has the same support and sign.
These results are stated as the following lemma.

%

\begin{lemma}
\label{Lemm_SuppConv}
Let $\{x^n\}$ be a sequence generated by IJT
algorithm. Assume that $0<\mu<\frac{1}{L}$,
then there exist a sufficiently large positive
integer $n^*$, an index set $I$ and a sign vector $S^*$ such that when $n > n^*$, it holds
\begin{enumerate}
\item[(a)]
$I^n=I$;

\item[(b)]
$Supp(x^*) = I, \forall x^* \in \cal{X}$;

\item[(c)]
$sign(x^n) = S^*$;

\item[(d)]
$sign(x^*) = S^*, \forall x^* \in \cal{X}$.
\end{enumerate}
\end{lemma}

The proof of this lemma is presented in Appendix C.
This lemma gives a possible way to construct a new sequence on a special subspace that has the same convergence behavior of $\{x^n\}$.
Thus, if we can prove the convergence of the new sequence, then the strong convergence of $\{x^n\}$ can naturally be claimed.
Specifically, such new sequence can be constructed as follows.
By Lemma {\ref{Lemm_SuppConv}}, there exists a sufficiently large integer $n^*>0$ such that when $n>n^*$,
\[
I^n = I \ \text{and} \ sign(x^n) = sign(x^*).
\]
Therefore, we can claim that $\{x^n\}$ converges to $x^*$ if the new sequence $\{x^{i+n^*}\}_{i\in \mathbf{N}}$ converges to $x^*$,
which is also equivalent to the convergence of the sequence $\{z^{i+n^*}\}_{i\in \mathbf{N}}$, i.e.,
\begin{equation}
z^{i+n^*} \rightarrow z^* \ \text{as}\  i\rightarrow \infty
\label{Seq1}
\end{equation}
with $z^{i+n^*} = P_I x^{i+n^*}$ and $z^* = P_I x^*$.
Let $\hat{z}^n = z^{n+n^*}$, then $\{\hat{z}^n\}$ has the same convergence behavior of $\{x^n\}$.


For any $\epsilon>0$, we define a one-dimensional real space
$$
\mathbf{R}_{\epsilon} = \mathbf{R}\setminus (-\epsilon, \epsilon).
$$
Particularly, let $\mathbf{R}_0 = \mathbf{R}\setminus \{0\}$.
Denote ${\cal{Z}}^*=P_I{\cal{X}}= \{P_Ix^*: x^* \in {\cal{X}}\}$.
We define two new functions $T: \mathbf{R}_{\eta_{\mu}/2}^K \rightarrow \mathbf{R}$ and $f: \mathbf{R}_{\eta_{\mu}/2}^K \rightarrow \mathbf{R}$ with
\begin{equation}
T(z) = T_{\lambda}(P_I^T z) \text{\ and\ } f(z) = F(P_I^T z),
\label{Def_T}
\end{equation}
for any $z\in \mathbf{R}_{\eta/2}^K$, respectively.
For any $z^* \in {\cal{Z}}^*$, it can be observed that $z^* \in \mathbf{R}_{\eta_{\mu}}^K$ by Lemma \ref{Lemm_JumpThresFun},
and $z^*$ is indeed a critical point of $T$ from Property \ref{Prop_OptCond}(a).
Moreover, we define a series of mappings $\phi_{1,m}: \mathbf{R}_0^m \rightarrow \mathbf{R}^m$
and $\phi_{2,m}: \mathbf{R}_0^m \rightarrow \mathbf{R}^{m\times m}$ as follows
\begin{align}
\label{phi1}
\phi_{1,m}(z) = (sign(z_1)\phi'(|z_1|),\cdots,sign(z_m)\phi'(|z_m|))^T, &
\end{align}
\begin{align}
\label{phi2}
\phi_{2,m}(z) = diag(\phi''(|z_1|),\cdots,\phi''(|z_m|)), &
\end{align}
$m=1,\cdots,N$,
where $diag(z)$ represents the diagonal matrix generated by $z$.
For brevity, we will denote $\phi_{1,m}$ and $\phi_{2,m}$ as $\phi_1$ and $\phi_2$ respectively when $m$ is fixed and there is no confusion.

By Properties {\ref{Prop_SuffDecrease}}-{\ref{Prop_OptCond}}, we can easily justify that $\{\hat{z}^n\}$ satisfies the following
so-called sufficient decrease, relative error and continuity conditions.

\begin{lemma}
\label{Lemm_NewSeq3Cond}
$\{\hat{z}^n\}$ satisfies the following conditions:
\begin{enumerate}
\item[(a)] (Sufficient decrease condition). For each $n\in \mathbf{N}$,
$$T(\hat{z}^{n+1}) \leq T(\hat{z}^n) - \frac{1}{2}(\frac{1}{\mu}-L) \|\hat{z}^{n+1}-\hat{z}^n\|_2^2.$$

\item[(b)] (Relative error condition). For each $n\in \mathbf{N}$,
$$\|\nabla T(\hat{z}^{n+1})\|_2 \leq (\frac{1}{\mu} + L) \|\hat{z}^{n+1} - \hat{z}^n\|_2.$$

\item[(c)] (Continuity condition). There exists a subsequence $\{\hat{z}^{n_j}\}_{j\in \mathbf{N}}$ and $z^*$ such that
$$
\hat{z}^{n_j} \rightarrow z^* \ \text{and}\ T(\hat{z}^{n_j}) \rightarrow T(z^*), \ \text{as}\ j\rightarrow \infty.
$$
\end{enumerate}
\end{lemma}

From this lemma, if $T$ further has the KL property at the limit point $z^*$, then according to Theorem 2.9 in {\cite{Attouch2013}},
$\{\hat{z}^n\}$ definitely converges to $z^*$.
Lemma {\ref{Lemm_NewSeq3Cond}}(a) and (c) are obvious by Properties {\ref{Prop_SuffDecrease}}-{\ref{Prop_SubseqConv}}, the specific form of $T$ and the construction of $\{\hat{z}^n\}$.
Lemma {\ref{Lemm_NewSeq3Cond}}(b) holds mainly due to Property {\ref{Prop_OptCond}}(b) and Assumptions \ref{AonF}-\ref{AonPhi}.
Specifically, by Property {\ref{Prop_OptCond}}(b), it can be easily checked that
\begin{align*}
\hat{z}^{n+1} + \lambda \mu\phi_1(\hat{z}^{n+1}) = \hat{z}^{n} -\mu \nabla f(\hat{z}^{n}),
\end{align*}
which implies
\begin{align*}
& \mu(\nabla f(\hat{z}^{n+1}) + \lambda \phi_1(\hat{z}^{n+1})) = \nonumber\\
&(\hat{z}^{n} - \hat{z}^{n+1}) + \mu(\nabla f(\hat{z}^{n+1})-\nabla f(\hat{z}^{n})).
\end{align*}
Thus,
\begin{align*}
\|\nabla T(\hat{z}^{n+1})\|_2 =\frac{1}{\mu} \|(\hat{z}^{n} - \hat{z}^{n+1}) + \mu(\nabla f(\hat{z}^{n+1})-\nabla f(\hat{z}^{n}))\|_2.
\end{align*}
By Assumption 1, $\nabla F$ is Lipschitz continuous with the Lipschitz constant $L$, then
\begin{align*}
&\|\nabla f(\hat{z}^{n+1})-\nabla f(\hat{z}^{n})\|_2 \\\nonumber
&= \|[\nabla F(P_I^T\hat{z}^{n+1})]_I - [\nabla F(P_I^T\hat{z}^{n})]_I\|_2 \nonumber\\
&\leq \|\nabla F(P_I^T\hat{z}^{n+1}) - \nabla F(P_I^T\hat{z}^{n})\|_2 \\\nonumber
&\leq L\|P_I^T\hat{z}^{n+1} - P_I^T\hat{z}^{n}\|_2
= L \|\hat{z}^{n+1} - \hat{z}^{n}\|_2.
\end{align*}
Therefore,
\begin{equation*}
\|\nabla T(\hat{z}^{n+1})\|_2 \leq (\frac{1}{\mu} + L) \|\hat{z}^{n+1} - \hat{z}^{n}\|_2.
\end{equation*}

By Lemma {\ref{Lemm_NewSeq3Cond}} and the construction form of $\{\hat{z}^n\}$, we can obtain the following convergence result of IJT algorithm.
\begin{thm}
\label{Thm_StrongConv}
Assume that $F$ and $\phi$ satisfy Assumptions 1 and 2, respectively.
Consider any sequence $\{x^n\}$ generated by IJT algorithm with a bounded initialization.
Suppose that $0<\mu<\frac{1}{L}$, then $\{x^n\}$ converges subsequentially to a set $\cal{X}$.
If further $T_{\lambda}$ satisfies the $I$-rKL property at some limit point $x^* \in {\cal{X}}$ with $I = Supp(x^*)$,
then the whole sequence $\{x^n\}$ indeed converges to $x^*$.
\end{thm}

The first part of this theorem states that the sequence $\{x^n\}$ converges subsequentially to a limit point set $\cal{X}$
as long as the step size parameter $\mu$ is sufficiently small.
The second part shows that the objective function further satisfies the introduced rKL property at some limit point $x^*$,
then the sequence $\{x^n\}$ converges to $x^*$.

Furthermore, combining Lemma {\ref{Lemm_SuffCond_rKL}} and Theorem {\ref{Thm_StrongConv}}, we can obtain the following corollary.

\begin{coro}
\label{Coro_StrongConv}
Assume that $F$ and $\phi$ satisfy Assumptions 1 and 2, respectively.
Consider any sequence $\{x^n\}$ generated by IJT algorithm with a bounded initialization.
Suppose that $0<\mu<\frac{1}{L}$, and
if further there exists a limit point $x^*$ such that $F$ is twice continuously differentiable at $x^*$ and $\nabla^2 T(P_Ix^*)$ is nonsingular,
then the whole sequence $\{x^n\}$ indeed converges to $x^*$.
\end{coro}

\subsection{Convergence To A Local Minimizer}

As shown in Corollary {\ref{Coro_StrongConv}}, if $\nabla^2 T(P_Ix^*)$ is nonsingular at some limit point $x^*$, then the sequence
generated by IJT algorithm converges to $x^*$, which is also a stationary point. In this subsection, we will justify that $x^*$ is also a
local minimizer of the optimization problem if $\nabla^2 T(P_Ix^*)$ is positive definite.

\begin{thm}
\label{Thm_LocalConv}
Suppose that $F$ and $\phi$ satisfy Assumptions 1
and 2, respectively. Assume that $0<\mu<\frac{1}{L}$, and the
sequence $\{x^n\}$ generated by IJT algorithm converges to
$x^*$. 
Then $x^*$
is a local minimizer of $T_{\lambda}$ provided that $F$ is twice continuously differentiable at $x^*$ and
$\nabla^2 T(P_Ix^*)$ is positive definite.
\end{thm}

The proof of this theorem is rather intuitive.
In the following, we will present some simple derivations.
By Property {\ref{Prop_OptCond}}(a) we have
\begin{equation}
[\nabla F(x^*)]_I + \lambda \phi_1(x_I^*) = 0.
\label{OptCond}
\end{equation}
This together with the condition of the theorem
$$
\nabla^2 T(P_Ix^*) = \nabla_{II}^2 F(x^*)+\lambda \phi_2(x_I^*) \succ 0
$$
imply that the second-order optimality conditions hold at $x^* = (x_I^*,0)$,
where $\nabla_{II}^2 F(x^*)={\frac{\partial^2 F(x)}{\partial x_I^2}}{\big|}_{x=x^*}.$
For sufficiently small vector $h$, we denote $x_h^* = (x_I^*+h_I,0)$.
It then follows
\begin{equation}
F(x_h^*)+\lambda \sum_{i\in I} \phi(|x_i^*+h_i|) \geq F(x^*) + \lambda \sum_{i\in I} \phi(|x_i^*|).
\label{SuppPart}
\end{equation}
Furthermore, by Assumption {\ref{AonPhi}}(c), it obviously holds that
$$
\phi(t)>(\|[\nabla F(x^*)]_{I^c}\|_{\infty}+2)t/\lambda,
$$
for sufficiently small $t>0$.
By this fact and the differentiability of $F$, one can observe that for sufficiently small $h$, there hold
\begin{align}
\label{ZeroPart}
&F(x^*+h)-F(x_h^*) + \lambda \sum_{i\in I^c} \phi(|h_i|)\nonumber\\
&=h_{I^c}^T [\nabla F(x^*)]_{I^c} + \lambda \sum_{i\in I^c} \phi(|h_i|)+ o(h_{I^c}) \nonumber\\
&\geq\sum_{i\in I^c} (\|[\nabla F(x^*)]_{I^c}\|_{\infty} - [\nabla  F(x^*)]_i+1)|h_i| \geq 0.
\end{align}
Summing up the above two inequalities (\ref{SuppPart})-(\ref{ZeroPart}),
one has that for all sufficiently small $h$,
\begin{equation}
T_{\lambda}(x^*+h) - T_{\lambda}(x^*) \geq 0,
\label{LocalMinimizer}
\end{equation}
and hence $x^*$ is a local minimizer.

Actually, we can observe that when $h \neq 0$,
then at least one of these two inequalities (\ref{SuppPart}) and (\ref{ZeroPart}) will hold strictly,
which implies that $x^*$ is a strictly local minimizer.

\subsection{Asymptotically Linear Convergence Rate}


In order to derive the rate of convergence of IJT algorithm, we
first show some observations
on $\nabla F$ and $\phi'$ in the neighborhood of $x^*$.
For any $0<\varepsilon<\eta_{\mu}$, we define a neighborhood of $x^*$ as follows
$$
{\cal{N}}(x^*,\varepsilon) = \{x\in \mathbf{R}^N: \|x_I - x_I^*\|_2< \varepsilon, x_{I^c} = 0 \}.
$$
If $F$ is twice continuously differentiable at $x^*$ and also $\lambda_{\min}(\nabla_{II}^2 F(x^*))>0$,
then for any $x\in {\cal{N}}(x^*,\varepsilon)$,
there exist two sufficiently small positive constants $c_F$ and $c_{\phi}$
(both $c_F$ and $c_{\phi}$ depending on $\varepsilon$ with $c_F \rightarrow 0$ and $c_{\phi} \rightarrow 0$ as $\varepsilon \rightarrow 0$)
such that
\begin{align}
\label{CondF}
&\langle [\nabla F(x)]_I - [\nabla F(x^*)]_I, x_I-x_I^* \rangle \\\nonumber
&\geq (\lambda_{\min}(\nabla_{II}^2 F(x^*)) - c_F) \|x_I-x_I^*\|_2^2,
\end{align}
and
\begin{align}
\label{Condphi}
&\langle \phi_1(x_I) - \phi_1(x_I^*), x_I-x_I^* \rangle \\\nonumber
&\geq (\phi''(e) - c_{\phi}) \|x_I-x_I^*\|_2^2,
\end{align}
where (\ref{Condphi}) holds for $\phi'$ being strictly convex on $(0,\infty)$, and thus $\phi''$ being nondecreasing on $(0,\infty)$,
consequently, $\min_{i\in I} \phi''(|x_i^*|) = \phi''(\min_{i\in I} |x_i^*|) $.
With the observations (\ref{CondF}) and (\ref{Condphi}), we obtain the following theorem.

\begin{thm}
\label{Thm_ConvRate1}
Suppose that $F$ and $\phi$ satisfy Assumptions 1
and 2, respectively. Assume that the sequence $\{x^n\}$ generated by
IJT algorithm converges to $x^*$.
Let
$e=\min_{i\in I} |x_i^*|$.
Moreover, if $F$ is twice continuously differentiable at $x^*$ and the following conditions
hold
\begin{enumerate}
\item[(a)]
$\lambda_{\min}(\nabla_{II}^2 F(x^*))>0$;

\item[(b)]
$0<\lambda < -\frac{\lambda_{\min}(\nabla_{II}^2 F(x^*))}{\phi''(e)},$

\item[(c)]
$0<\mu<\min \{\frac{2(\lambda_{\min}(\nabla_{II}^2 F(x^*)) + \lambda \phi''(e))}{L^2 - (\lambda \phi''(e))^2}, \frac{1}{L}\}$,
\end{enumerate}
then there exists a sufficiently large positive integer $n_{0}$ and
a constant $\rho^{\ast}\in(0,1)$ such that when $n>n_{0},$
\[
\Vert x^{n+1}-x^{\ast}\Vert_{2}\leq\rho^{\ast}\Vert x^{n}-x^{\ast}%
\Vert_{2},
\]
and
\[
\|x^{n+1}-x^*\|_2 \leq \frac{\rho^*}{1-\rho^*} \|x^{n+1} - x^n\|_2.
\]
\end{thm}

The proof of Theorem {\ref{Thm_ConvRate1}} is presented in Appendix D.
This theorem states that IJT algorithm has asymptotically linear convergence rate under certain conditions.
Let $z^* = P_I x^*$.
Conditions (a) and (b) in this theorem imply that the Hessian of $T$ at $z^*$, $\nabla^2 T(z^*)$ is strongly positive definite,
since $\lambda_{\min}(\nabla^2 T(z^*)) = \lambda_{\min}(\nabla^2 f(z^*) + \lambda \phi_2(z^*))
\geq \lambda_{\min}(\nabla^2 f(z^*)) + \lambda \cdot \lambda_{\min}(\phi_2(z^*)) =  \lambda_{\min}(\nabla^2 f(z^*)) + \lambda \phi''(e)>0.$
Thus, $T$ is locally strongly convex at $z^*$.
Theorem {\ref{Thm_ConvRate1}} actually implies that the auxiliary sequence $\{\hat{z}^n\}$ converges linearly
if $T$ is strongly convex at $z^*$ and the step size parameter $\mu$ is sufficiently small.
As shown by this theorem, if we can fortunately obtain a sufficiently good initialization,
then IJT algorithm may converge fast with a linear rate.
On the other hand, Theorem {\ref{Thm_ConvRate1}} also provides a posteriori computable error estimation of the algorithm,
which can be used to design an efficient terminal rule of IJT algorithm.

%
%


It can be observed that the conditions of Theorem {\ref{Thm_ConvRate1}} are slightly stricter than those of Corollary {\ref{Coro_StrongConv}},
and thus, $x^*$ is also a local minimizer under the conditions of Theorem {\ref{Thm_ConvRate1}}.
In the following, we will show that the condition on $\mu$ in Theorem {\ref{Thm_ConvRate1}} can be extended to $0<\mu<1/L$ if we add some additional assumptions on the higher order differentiability of $\phi$ in the neighborhood of the local minimizer $x^*$.
We state this as the following theorem.

\begin{thm}
\label{Thm_ConvRate2}
Assume that $0<\mu< \frac{1}{L}$. Let $\{x^n\}$ be
a sequence generated by IJT algorithm and converge to $x^*$. Let
$e=\min_{i\in I} |x_i^*|$.
Moreover, if $F$ is twice continuously differentiable at $x^*$ and the following conditions hold
\begin{enumerate}
\item[(a)]
$\lambda_{\min}(\nabla_{II}^2 F(x^*))>0$,

\item[(b)]
$0<\lambda < -\frac{\lambda_{\min}(\nabla_{II}^2 F(x^*))}{\phi''(e)},$

\item[(c)]
for any sufficiently small $0<\varepsilon<\eta_{\mu}$, the derivative of $\phi''$, $\phi'''$ is well-defined, bounded and nonzero on the set $\cup_{i\in I} B(x_i^*,\varepsilon)$, where $B(x_i^*,\varepsilon) := (x_i^*-\varepsilon, x_i^*+\varepsilon)$,
\end{enumerate}
then there exists a sufficiently large positive integer $n_0>0$ and
a constant $\rho\in(0,1)$ such that when $n>n_0,$
\[
\Vert x^{n+1}-x^{\ast}\Vert_{2}\leq\rho\Vert x^{n}-x^{\ast}%
\Vert_{2},
\]
and
\[
\|x^{n+1}-x^*\|_2 \leq \frac{\rho}{1-\rho}\|x^{n+1}-x^n\|_2.
\]
\end{thm}

The proof of this theorem is given in Appendix E.
Note that the condition (c) can be easily satisfied if the penalty $\phi$ has the continuous third-order derivative on $(0,\infty)$.
In the next section, we will show that the $l_q$-norm ($0<q<1$) is one of the most typical subclass of these non-convex penalties that satisfy the condition (c) in Theorem {\ref{Thm_ConvRate2}}.

\section{Application to $l_q$ Regularization ($0<q<1$)}

In this section, we apply the established theoretical results to a typical case, $l_q$ regularization with $0<q<1$.

Mathematically, $l_q $ $(0<q<1)$ regularization can be formulated as follows
\begin{equation*}
\min_{x\in \mathbf{R}^N} \left\{T_{\lambda}(x) = \frac{1}{2}\|Ax-y\|_2^2 + \lambda \|x\|_q^q \right\},
\end{equation*}
where $A\in \mathbf{R}^{M\times N}$ (commonly, $M<N$) is usually called the sensing matrix,
$y\in \mathbf{R}^M$ is called the measurement vector, $x$ is commonly assumed to be sparse, i.e., $\|x\|_0 \ll N$,
and $\|x\|_q^q = \sum_{i=1}^N |x_i|^q$.
Thus, in such special case, $F(x) = \frac{1}{2}\|Ax-y\|_2^2$ and $\Phi(x) = \|x\|_q^q$ with $\phi(x) = x^q$ defined on $(0,\infty)$.
In {\cite{BrediesNonconvex}}, Bredies and Lorenz demonstrated that the one-dimensional proximity operator $prox_{\mu, \lambda |\cdot|^q}$ of $l_q$-norm can be expressed as
\begin{equation}
prox_{\mu, \lambda |\cdot|^q} (z)=\left\{
\begin{array}
[c]{ll}%
(\cdot + \lambda \mu q sign(\cdot) |\cdot|^{q-1})^{-1}(z), & |z|\geq \tau_{\mu,q}\\
0, &  |z|\leq \tau_{\mu,q}
\end{array}
\right.
\label{ProxMapExpLq}
\end{equation}
for any $z \in \mathbf{R}$ with
\begin{equation}
\tau_{\mu,q} = \frac{2-q}{2-2q}(2\lambda\mu(1-q))^{\frac{1}{2-q}},
\label{ThreshValuexLq}
\end{equation}
\begin{equation}
\eta_{\mu,q} = (2\lambda\mu (1-q))^{\frac{1}{2-q}},
\label{ThreshValueyLq}
\end{equation}
and the range of $prox_{\mu, \lambda |\cdot|^q}$ is $\{0\}\cup [\eta_{\mu,q},\infty)$.
Furthermore, for some special $q$ (say, $q=1/2, 2/3$),
the corresponding proximity operators can be expressed analytically {\cite{L1/2TNN}}, {\cite{L2/3Cao2013}}.

According to {\cite{Attouch2013}} (See Example 5.4, page 122), the function $T_{\lambda}(x) = \frac{1}{2}\|Ax-y\|_2^2 + \lambda \|x\|_q^q$
is a KL function and obviously satisfies the rKL propety at any limit point.
By applying Theorem {\ref{Thm_StrongConv}} to the $l_q$ regularization, we can obtain the following corollary directly.

\begin{coro}
\label{Coro_StrongConvLq}
Let $\{x^n\}$ be a sequence generated by IJT algorithm for $l_q$ regularization with $q\in (0,1)$.
Assume that $0<\mu<\frac{1}{\|A\|_2^2}$, then $\{x^n\}$ converges to a stationary point of $l_q$ regularization.
\end{coro}

In {\cite{Attouch2013}}, Attouch et al. showed the convergence of the inexact forward-backward splitting algorithm for $l_q$ regularization (See Theorem 5.1, page 118) under exactly the same condition of Corollary {\ref{Coro_StrongConvLq}} .
Furthermore, it is easy to check that $F(x) = \frac{1}{2}\|Ax-y\|_2^2$ and $\phi(z) = z^q$ satisfy Assumptions 1 and 2, respectively.
In addition, $\phi(z) = z^q$ also satisfies the condition (c) in Theorem {\ref{Thm_ConvRate2}} naturally.
Therefore, as a direct corollary of Theorem {\ref{Thm_ConvRate2}}, we show the asymptotically linear convergence rate of IJT algorithm for $l_q$ regularization as follows.

\begin{coro}
\label{Coro_ConvRateLq}
Assume that $0<\mu<\|A\|_2^{-2}$.
Let $\{x^n\}$ be a sequence generated by IJT algorithm for $l_q$ ($0<q<1$) regularization and converge to $x^*$.
Let $I=Supp(x^*)$ and $e=\min_{i\in I} |x_i^*|$.
Moreover, if the following conditions hold:
\begin{enumerate}
\item[(a)]
$\lambda_{\min}(A_I^TA_I)>0$,

\item[(b)]
$0<\lambda < \frac{\lambda_{\min}(A_I^TA_I)e^{2-q}}{q(1-q)},$
\end{enumerate}
then there exists a sufficiently large positive integer $n_0$ and
a constant $\rho\in(0,1)$ such that when $n>n_0,$
\[
\Vert x^{n+1}-x^{\ast}\Vert_{2}\leq\rho\Vert x^{n}-x^{\ast}%
\Vert_{2},
\]
and
\[
\|x^{n+1} - x^*\|_2 \leq \frac{\rho}{1-\rho} \|x^{n+1}-x^n\|_2.
\]
In addition, $x^*$ is also a local minimizer of $l_q$ regularization.
\end{coro}

The condition (b) in Corollary {\ref{Coro_ConvRateLq}} means that the regularization parameter should be sufficiently small to guarantee that the limit point is a local minimizer.
Instead of adding the assumption on the regularization parameter $\lambda$,
we give another sufficient condition characterized by the matrix $A$.
Such condition is mainly derived via taking advantage of the specific form of the threshold value (\ref{ThreshValueyLq}).
More specifically, by (\ref{ThreshValueyLq}), it holds
\begin{equation}
e\geq \eta_{\mu,q} = (2\lambda\mu (1-q))^{\frac{1}{2-q}}.
\end{equation}
Then if $\frac{\lambda_{\min}(A_I^TA_I)}{\|A\|_2^2} > \frac{q}{2}$ and
\begin{equation}
\frac{q}{2\lambda_{\min}(A_I^TA_I)} < \mu < \frac{1}{\|A\|_2^2},
\label{SuffCond2}
\end{equation}
the conditions in Corollary {\ref{Coro_ConvRateLq}} hold naturally. Therefore, we can obtain the following theorem on the asymptotically linear convergence rate of IJT algorithm applied to $l_q$ regularization.

\begin{thm}
\label{Thm_ConvRateLq_mu}
Assume that $0<\mu<\|A\|_2^{-2}$.
Let $\{x^n\}$ be a sequence generated by IJT algorithm for $l_q$ ($0<q<1$) regularization and converge to $x^*$.
Let $I=Supp(x^*)$.
Moreover, if the following conditions hold:
\begin{enumerate}
\item[(a)]
$\frac{\lambda_{\min}(A_I^TA_I)}{\|A\|_2^2} > \frac{q}{2}$,

\item[(b)]
$\frac{q}{2\lambda_{\min}(A_I^TA_I)} < \mu < \frac{1}{\|A\|_2^2},$
\end{enumerate}
then there exists a sufficiently large positive integer $n_0$ and
a constant $\rho\in(0,1)$ such that when $n>n_0,$
\[
\Vert x^{n+1}-x^{\ast}\Vert_{2}\leq\rho\Vert x^{n}-x^{\ast}%
\Vert_{2},
\]
and
\[
\|x^{n+1} - x^*\|_2 \leq \frac{\rho}{1-\rho} \|x^{n+1}-x^n\|_2.
\]
In addition, $x^*$ is also a local minimizer of $l_q$ regularization.
\end{thm}

From Theorem {\ref{Thm_ConvRateLq_mu}}, it means that if the matrix $A$ satisfies a certain concentration property and the step size $\mu$ is chosen appropriately,
then IJT algorithm can converge to a local minimizer at an asymptotically linear rate.
Note that the condition (a) in Theorem {\ref{Thm_ConvRateLq_mu}} implies $\frac{q}{2\lambda_{\min}(A_I^TA_I)}<\frac{1}{\|A\|_2^2}$ naturally.
Thus, the condition (b) of Theorem {\ref{Thm_ConvRateLq_mu}} is a natural and reachable condition and,
furthermore, whenever this condition is satisfied, the sequence $\{x^n\}$ is indeed convergent by Corollary {\ref{Coro_StrongConvLq}}.
This shows that only the condition (a) is essential in Theorem {\ref{Thm_ConvRateLq_mu}}.
We notice that the condition (a) is a concentration condition on eigenvalues of the submatrix $A_I^TA_I$,
and, in particular, it implies
$$
\lambda_{\min}(A_I^TA_I) > q\lambda_{\max}(A_I^TA_I)/2,
$$
or equivalently
\begin{equation}
Cond(A_I^TA_I) : = \frac{\lambda_{\max}(A_I^TA_I)}{\lambda_{\min}(A_I^TA_I)} < \frac{2}{q},
\label{CondNum}
\end{equation}
where $Cond(A_I^TA_I)$ is the condition number of $A_I^TA_I$.
(\ref{CondNum}) thus shows that the submatrix $A_I^TA_I$ is well-conditioned with the condition number lower than $2/q$.

In recent years, a property called the restricted isometry property (RIP) of a
matrix $A$ was introduced to characterize the concentration degree of the
eigenvalues of its submatrix with $k$ columns {\cite{Candes05RIP}}. A
matrix $A$ is said to be of the $k$-order RIP (denoted then by $\delta_{k}%
$-RIP) if there exists a $\delta_{k}\in(0,1)$ such that
\begin{equation}
(1-\delta_{k})\Vert x\Vert_{2}^{2}\leq\Vert Ax\Vert_{2}^{2}\leq(1+\delta
_{k})\Vert x\Vert_{2}^{2},~\forall\Vert x\Vert_{0}\leq k. \label{RIP}%
\end{equation}
In other words, the RIP ensures that all submatrices of $A$ with $k$ columns
are close to an isometry, and therefore distance-preserving. Let $K=\Vert
x^{\ast}\Vert_{0}$. It can be seen from (\ref{RIP}) that if $A$ possesses
$\delta_{K}$-RIP with $\delta_{K}<\frac{2-q}{2+q}$, then
\[
Cond(A_{I}^{T}A_{I})\leq\frac{1+\delta_{K}}{1-\delta_{K}}<\frac{2}{q}.
\]
Thus, we can claim that when $A$ satisfies a certain RIP, the condition (a) in Theorem {\ref{Thm_ConvRateLq_mu}} can be satisfied. In particular, we have the
following proposition.

\begin{prop}
\label{Proposition_RIP}
Assume that $K<N/2$ and $A$ satisfies $\delta_{K}$-RIP
with $\delta_{K}<\frac{2-q}{2+2qN/K}$ or $\delta_{2K}$-RIP with $\delta
_{2K}<\frac{2-q}{2+qN/K}$, then the condition (a) in Theorem {\ref{Thm_ConvRateLq_mu}} holds.
\end{prop}

This can be directly checked by the facts that $\lambda_{\min}(A_{I}^{T}%
A_{I})\geq 1-\delta_{K}$, $\lambda_{\min}(A_{I}^{T}A_{I})\geq 1-\delta_{2K}$,
$\lambda_{\max}(A^{T}A)\leq 1+\delta_{N}$, $\delta_{N}\leq\frac{2N}{K}\delta
_{K}$ and $\delta_{N}\leq\frac{N}{K}\delta_{2K}$ (c.f. Proposition 1 in
{\cite{Foucart2010}}).

From Proposition {\ref{Proposition_RIP}}, we can see, for instance, when $q=1/2, K/N=1/3$ and $A$ satisfies
$\delta_{K}$-RIP with $\delta_{K}<3/10$ or $\delta_{2K}$-RIP with $\delta
_{2K}<3/7$, the condition (a) in Theorem {\ref{Thm_ConvRateLq_mu}} is satisfied, and therefore, by Theorem  {\ref{Thm_ConvRateLq_mu}}, IJT algorithm converges to a local minimizer of the $l_q$
regularization at an asymptotically linear rate. It is noted that in the condition of Proposition {\ref{Proposition_RIP}}, we always
have $\delta_{k}<\frac{2-q}{2+4q}$ and $\delta_{2k}<\frac{2-q}{2+2q}.$

\begin{remark}
\label{RHalf}
In a recent paper {\cite{ZengHalfConv2013}}, Zeng et al. have justified the convergence of a specific iterative thresholding algorithm called the iterative \textit{half} thresholding algorithm for $l_{1/2}$ regularization.
It can be observed that the convergence results of the iterative \textit{half} thresholding algorithm obtained in {\cite{ZengHalfConv2013}} is just a special case of the results presented in this section.
\end{remark}

\begin{remark}
\label{RHard1}
Recently, Lu {\cite{LU2014IHT}} proposed an iterative \textit{hard} thresholding method
and its variant for solving $l_0$ regularization over a conic constraint,
and established its convergence as well as the iteration complexity.
Although the $l_0$-norm does not satisfies Assumption {\ref{AonPhi}},
it can be observed that the finite support and sign convergence property (i.e., Lemma {\ref{Lemm_SuppConv}})
holds naturally for \textit{hard} algorithm due to
the \textit{hard} thresholding function possesses the similar discontinuity of the jumping thresholding function.
Furthermore, once the support of the sequence converges,
the iterative form of \textit{hard} algorithm is equal to the simple Landweber iteration,
and thus the convergence and asymptotically linear convergence rate of \textit{hard} algorithm can be directly claimed.
\end{remark}

\section{Related Work}

Recently, Attouch et al. {\cite{Attouch2013}} have justified the convergence of a family of descent methods
by assuming the objective function has the KL property {\cite{KLBolte2006}}, {\cite{KLBolte2007}},
and also the generated sequence satisfies the sufficient decrease property, relative error condition and continuity condition (Sec. 2.3 in {\cite{Attouch2013}}).
%
%
%
%
%
%
%
Instead of the well-known KL inequality condition, we introduce a weaker condition called the rKL property to check the convergence of IJT algorithm.
Besides the strong convergence, we also justify the asymptotically linear convergence rate of IJT algorithm under certain second-order conditions.
Compared with the
other algorithms including HQ {\cite{HQCon2006}}, FOCUSS
{\cite{FOCUSS}}, IRL1 {\cite{IRL1Chen}} and DC programming
{\cite{DCProgGasso2009}} algorithms, we derive a sufficient
condition instead of the direct assumption that the accumulation points are isolated,
for the convergence of IJT algorithm.
Furthermore, the convergence speed of IJT algorihtm is also
demonstrated in this paper.

Besides the aforementioned non-convex algorithms, there are some other related algorithms.
In the following, we will compare the obtained theoretical results of IJT algorithm with those of these algorithms.
The first class of closely related algorithms are the iterative shrinkage and thresholding (IST) algorithms,
which mainly refer to two generic algorithms and some specific algorithms.
The first generic algorithm related to IJT algorithm is
the generalized gradient projection (called GGP for short) algorithm {\cite{GGPYinTR2007}}, {\cite{BrediesNonconvex}}.
In {\cite{GGPYinTR2007}}, the GGP algorithm was proposed for the $l_1$ regularization problem.
In such a convex setting, the finite support convergence and eventually linear convergence rate was given in {\cite{GGPYinTR2007}}.
In {\cite{BrediesNonconvex}}, Bredies and Lorenz extended the GGP algorithm to solve the following general non-convex optimization model in the infinite-dimensional Hilbert space
\begin{equation}
\min_{x\in \mathbf{X}}\left\{F(x)+\lambda\Phi(x)\right\},
\label{NonconvexOpt}%
\end{equation}
where $\mathbf{X}$ is an infinite-dimensional Hilbert space,
$F: \mathbf{X} \rightarrow [ 0,\infty)$ is assumed to be a proper lower-semicontinuous
function with Lipschitz continuous gradient $\nabla F(x)$, and $\Phi: \mathbf{X} \rightarrow [ 0,\infty)$ is weakly lower-semicontinuous (possibly non-smooth and non-convex).
Furthermore, the iterative form of the GGP algorithm is specified as
$$
x^{n+1} \in Prox_{\mu,\lambda\Phi}(x^n-\mu \nabla F(x^n)),
$$
where $Prox_{\mu,\lambda\Phi}$ represents the proximity operator of $\Phi$ as defined in (\ref{ProxOper}).
It can be observed that IJT algorithm is a special case of GGP algorithm when applied to a separable $\Phi$ in the finite-dimensional real space.
Nevertheless, it was only
justified that GGP algorithm can converge subsequentially to a
stationary point {\cite{BrediesNonconvex}} (that is, there is a subsequence
that converges to a stationary point).
However, as a specific case of GGP
algorithm, we have justified that IJT algorithm can assuredly
converge to a local minimizer at an asymptotically linear convergence rate under certain conditions.

Another closely related generic algorithm is the general iterative
shrinkage and thresholding (GIST) algorithm suggested in {\cite{YeNonconvex}}.
The  GIST algorithm is proposed for the following general non-convex regularized optimization problem
\begin{equation}
\min_{x\in \mathbf{R}^N} \{F(x) + \lambda R(x)\},
\label{GISTAOpt}
\end{equation}
where $F$ is assumed to be continuously differentiable with Lipschitz continuous derivative,
and $R(x)$ is a continuous function and can be rewritten as the difference of two different convex functions.
As compared with Assumption {\ref{AonPhi}}, we can find that the optimization model considered in this paper is distinguished from the model (\ref{GISTAOpt}) studied in {\cite{YeNonconvex}}.
Moreover, only the subsequential convergence of the GIST algorithm can be justified in {\cite{YeNonconvex}},
while the convergence of the whole sequence and further the asymptotically linear convergence rate of IJT algorithm are demonstrated in this paper.

Besides these two generic algorithms, there are some other specific iterative thresholding algorithms related to IJT algorithm.
Among them, the \textit{hard} algorithm and the
\textit{soft} algorithm are two
representatives, which respectively solves the $l_{0}$ regularization and
$l_{1}$ regularization {\cite{Blumensath08}}, {\cite{Daubechies2004Soft}}.
It was demonstrated in {\cite{Blumensath08}}, {\cite{Daubechies2004Soft}} that when $\mu=1$
both $hard$ and $soft$ algorithms can converge to a stationary point whenever
$\Vert A\Vert_{2}<1$. These classical convergence results can be generalized
when a step size parameter $\mu$ is incorporated with the IST procedures, and in this case,
the convergence condition becomes
\begin{equation}
0<\mu<\Vert A\Vert_{2}^{-2}. \label{4.1}%
\end{equation}
It can be seen from Corollary {\ref{Coro_StrongConvLq}} that (\ref{4.1}) is the exact condition of the convergence of IJT algorithm when applied to the $l_q$ regularization with $0<q<1$,
which then supports that the classical convergence results of IST has been
extended to the non-convex $l_{q}$ ($0<q<1$) regularization case.
Furthermore, it was shown in {\cite{SoftBredies08}} that when the measurement matrix $A$
satisfies the so-called finite basis injective (FBI) property and the
stationary point possesses a strict sparsity pattern, the \textit{soft}
algorithm can converge to a global minimizer of $l_{1}$ regularization with a
linear convergence rate. Such result is
not surprising because of the convexity of $l_{1}$ regularization. As for
convergence speed of the \textit{hard} algorithm, it was demonstrated in {\cite{Blumensath08}} that
under the condition $\mu=1$ and $\Vert A\Vert_{2}<1$,
\textit{hard} algorithm will
converge to a local minimizer with an asymptotically linear convergence rate.
However, as algorithms for solving non-convex models,
Corollary {\ref{Coro_ConvRateLq}} and Theorem {\ref{Thm_ConvRateLq_mu}} reveal that IJT algorithm shares the same asymptotic
convergence speed with \textit{hard} algorithm.

\section{Numerical Experiments}

We conduct a set of numerical experiments in this section to substantiate the
validity of the theoretical analysis on the convergence of IJT algorithm.
While the effectiveness of IJT algorithm applied to large-scale applications such as the synthetic aperture radar (SAR) imaging and image processing can be referred to {\cite{SARZeng2012}} and {\cite{L2/3Cao2013}}.
(The corresponding matlab code of IJT algorithm can be referred to \url{https://github.com/JinshanZeng/IJT_Alg}.)

\subsection{Convergence Rate Justification}

\begin{figure}[!t]
\begin{minipage}[b]{.49\linewidth}
\centering
\includegraphics*[scale=0.33]{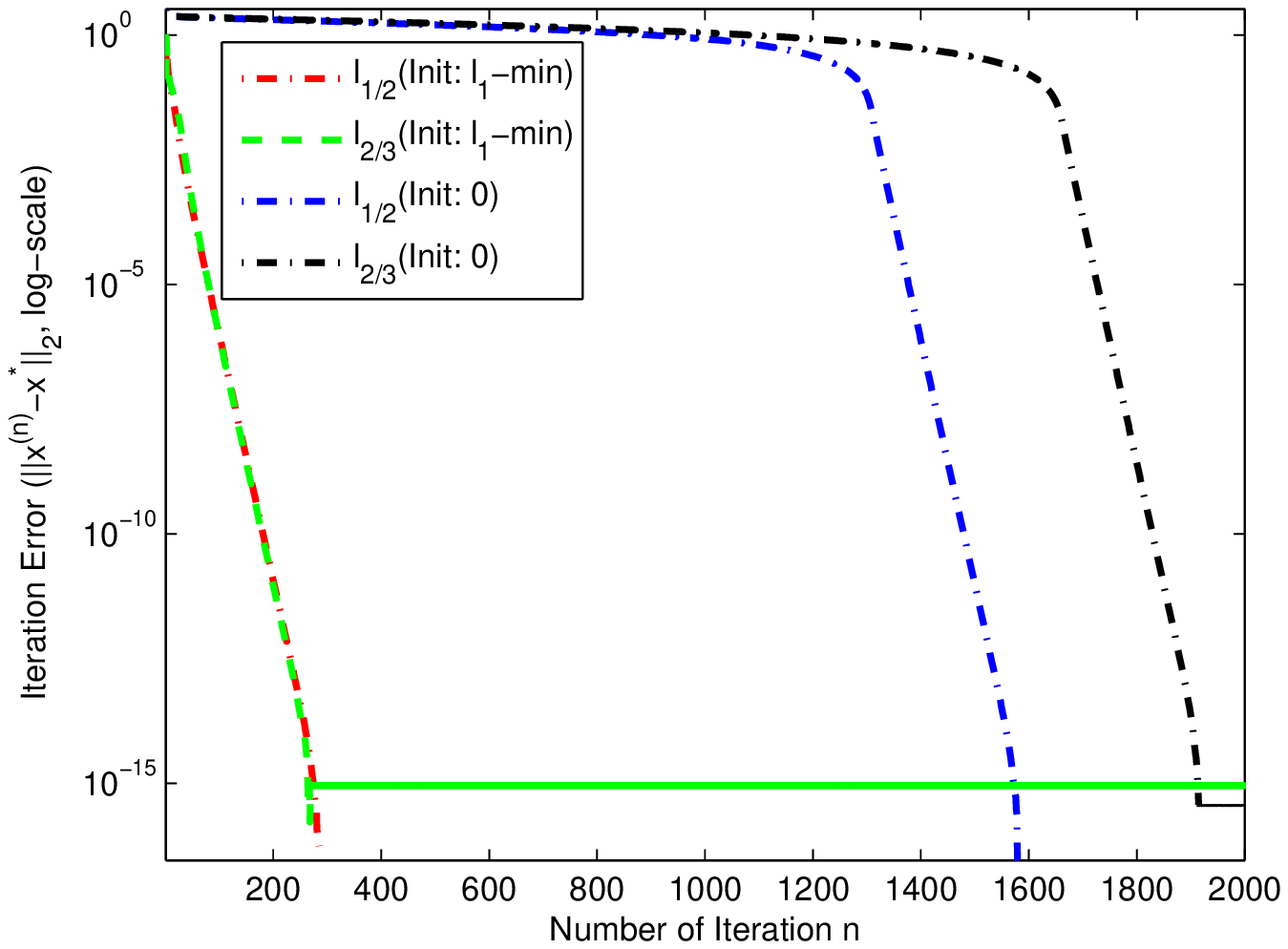}
\centerline{{\small (a) Iteration error}}
\end{minipage}
\hfill
\begin{minipage}[b]{.49\linewidth}
\centering
\includegraphics*[scale=0.33]{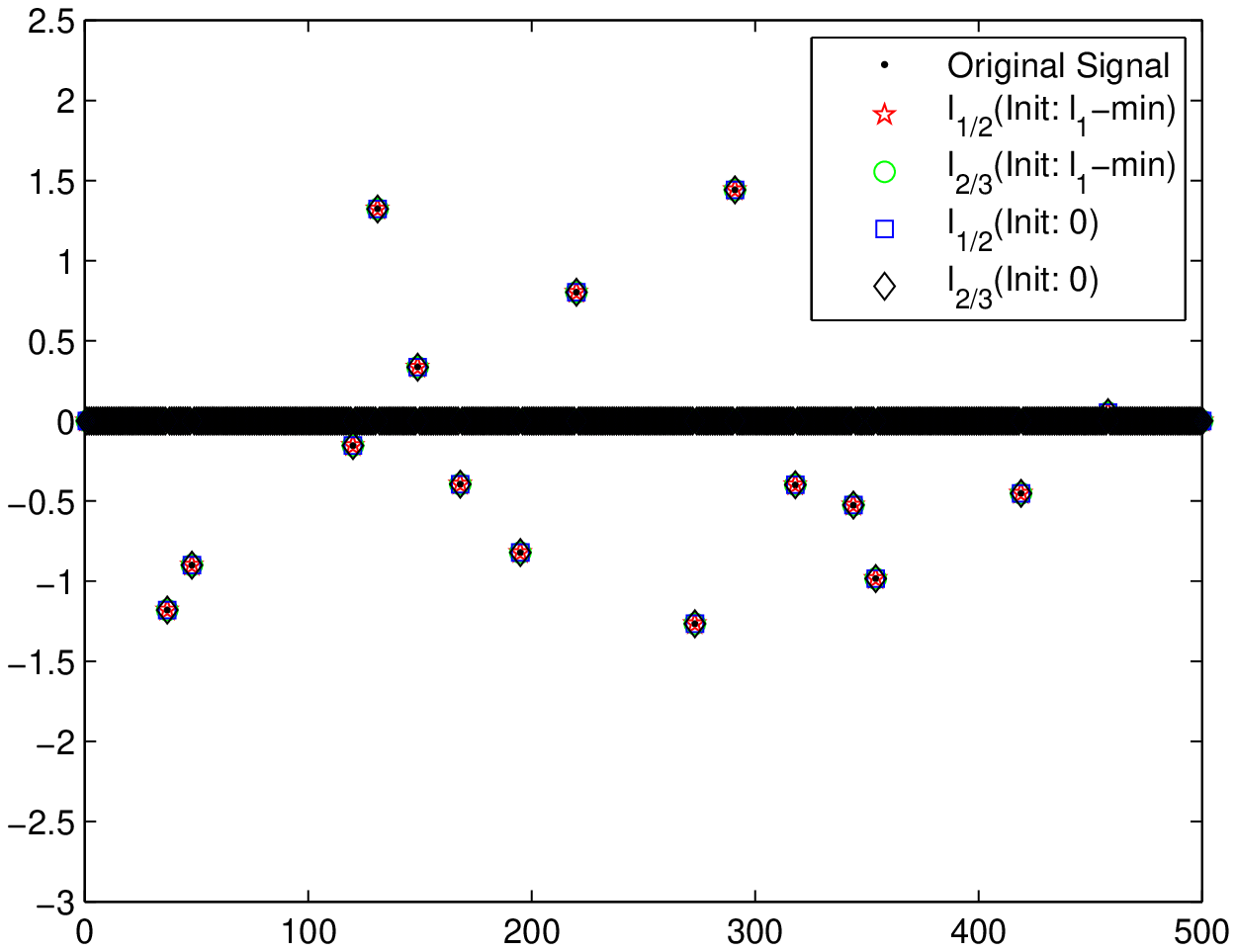}
\centerline{{\small (b) Recovery signal}}
\end{minipage}
\hfill
\caption{ Experiment for asymptotically linear convergence rate.
(a) The trend of iteration error, i.e., $\|x^{(n)}-x^*\|_2$.
(b) Recovery signal.
The labels ``$l_{1/2}$ (Init: $l_1$-min)'' and ``$l_{2/3}$ (Init: $l_1$-min)'' represent the cases of
$\phi(|z|)=|z|^{1/2}$ and $\phi(|z|) = |z|^{2/3}$ with the solution of the $l_1$-minimization problem as the initial guess, respectively.
The labels ``$l_{1/2}$ (Init: 0)'' and ``$l_{2/3}$ (Init: 0)'' represent the cases of
$\phi(|z|)=|z|^{1/2}$ and $\phi(|z|)=|z|^{2/3}$ with 0 as the initial guess, respectively.
The Recovery MSEs of the four cases, that is, $l_{1/2}$ (Init: $l_1$-min), $l_{2/3}$ (Init: $l_1$-min), $l_{1/2}$ (Init: 0) and $l_{2/3}$ (Init: 0)  are $3.06\times 10^{-6}$, $3.36\times 10^{-6}$, $3.24\times 10^{-6}$ and $3.67\times 10^{-6}$, respectively.
}
\label{Fig_AsympConvRate}
\end{figure}

We start with an experiment to confirm the linear rate of asymptotic
convergence. For this purpose, given a sparse signal $x$ with dimension
$N=500$ and sparsity $k=15,$ shown as in Fig. {\ref{Fig_AsympConvRate}}(b), we considered the signal
recovery problem through observation $y=Ax,$ where the measurement matrix $A$
is of dimension $M\times N=250\times500$ with Gaussian $\mathcal{N}(0,1/250)$
i.i.d. entries. Such measurement matrix is known to satisfy (with high
probability) the RIP with optimal bounds \cite{RudelsonRIP}, \cite{BaraniukRIP}.
We then applied IJT algorithm to the problem with two different non-convex penalties, that is, $\phi(|z|) = |z|^{1/2}, |z|^{2/3}$.
In both cases, the jumping thresholding operators can be analytically expressed as shown in {\cite{L1/2TNN}} and {\cite{L2/3Cao2013}}, respectively, and thus the corresponding IJT algorithms can be efficiently implemented.
In both cases, we took $\lambda=0.001$ and $\mu=0.99\Vert A\Vert_{2}^{-2}$.
Moreover, we considered two different initial guesses including 0 and the solution of the $l_1$-minimization problem
to justify the effect on the convergence speed.
The experiment results are reported in Fig. {\ref{Fig_AsympConvRate}}.

It can be seen from Fig. {\ref{Fig_AsympConvRate}}(a) how the iteration error ($\Vert
x^{(n)}-x^{\ast}\Vert_{2})$ varies.
More specifically, when 0 was taken as the initial guess, after approximately $1300$ and $1700$ iterations,
IJT algorithm converges to a stationary point with a linear decay rate for both penalties $\phi(|z|) = |z|^{1/2}$ and $\phi(|z|) = |z|^{2/3}$, as shown by the blue and black lines in Fig. {\ref{Fig_AsympConvRate}}(a), respectively.
While from the red and green lines in Fig. {\ref{Fig_AsympConvRate}}(a),
if we took the solution of the $l_1$-minimization problem as the initialization,
the IJT algorithm converges to a stationary point with a linear convergence rate starting from almost the first iteration for both penalties.
This indicates that the solution of the $l_1$-minimization problem is a good initialization, which is sufficiently close to the stationary point.
Moreover, Fig. {\ref{Fig_AsympConvRate}}(b) shows that the original sparse signal has been recovered by IJT algorithm with very high accuracy.
This experiment clearly justifies the convergence properties of IJT algorithm we have verified,
particularly the expected asymptotically linear convergence rate of IJT algorithm is substantiated.

\subsection{On effect of $\mu$}

\begin{figure}[!t]
\begin{minipage}[b]{.32\linewidth}
\centering
\includegraphics*[scale=0.21]{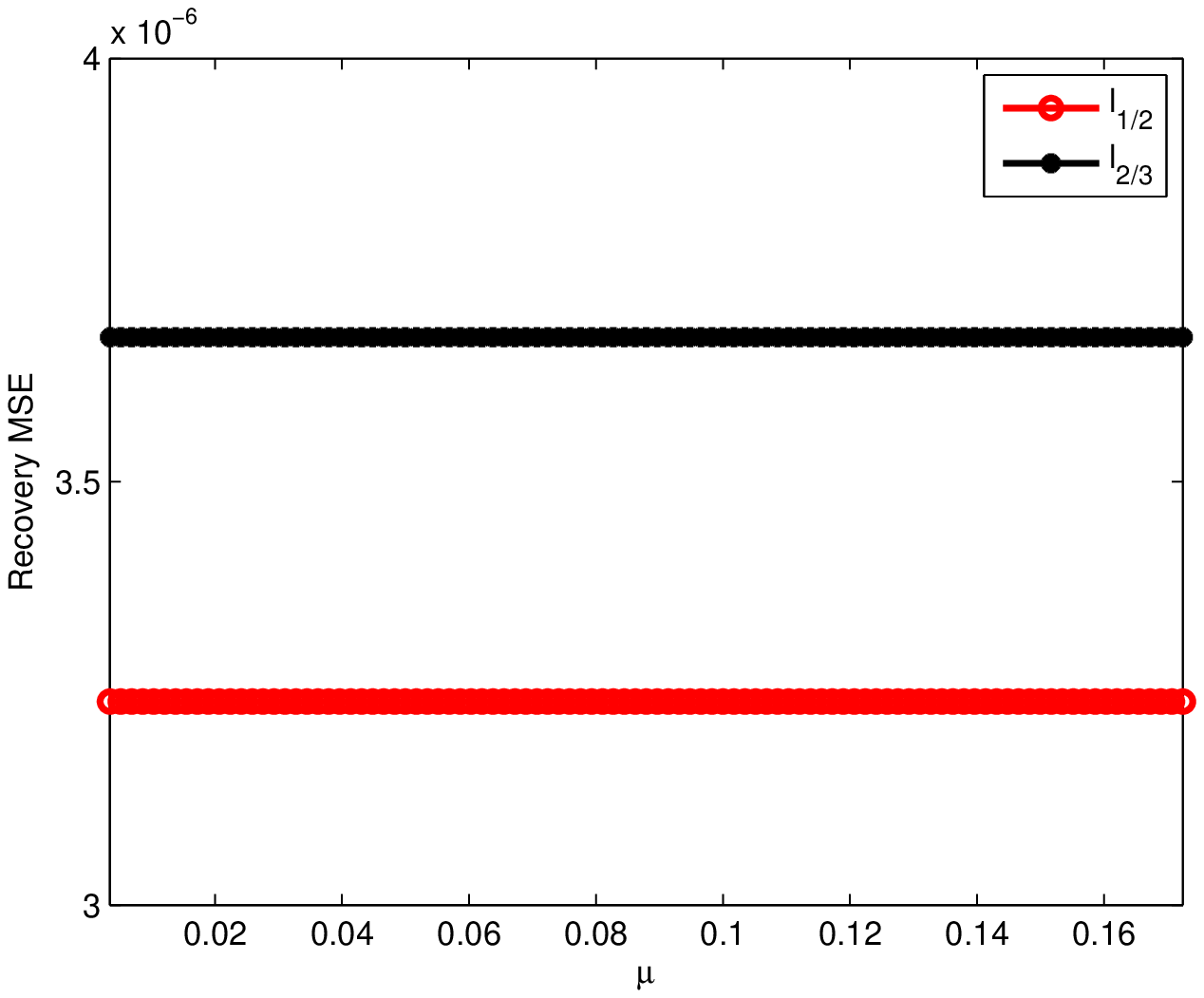}
\centerline{{\small (a) Recovery error}}
\end{minipage}
\begin{minipage}[b]{.32\linewidth}
\centering
\includegraphics*[scale=0.21]{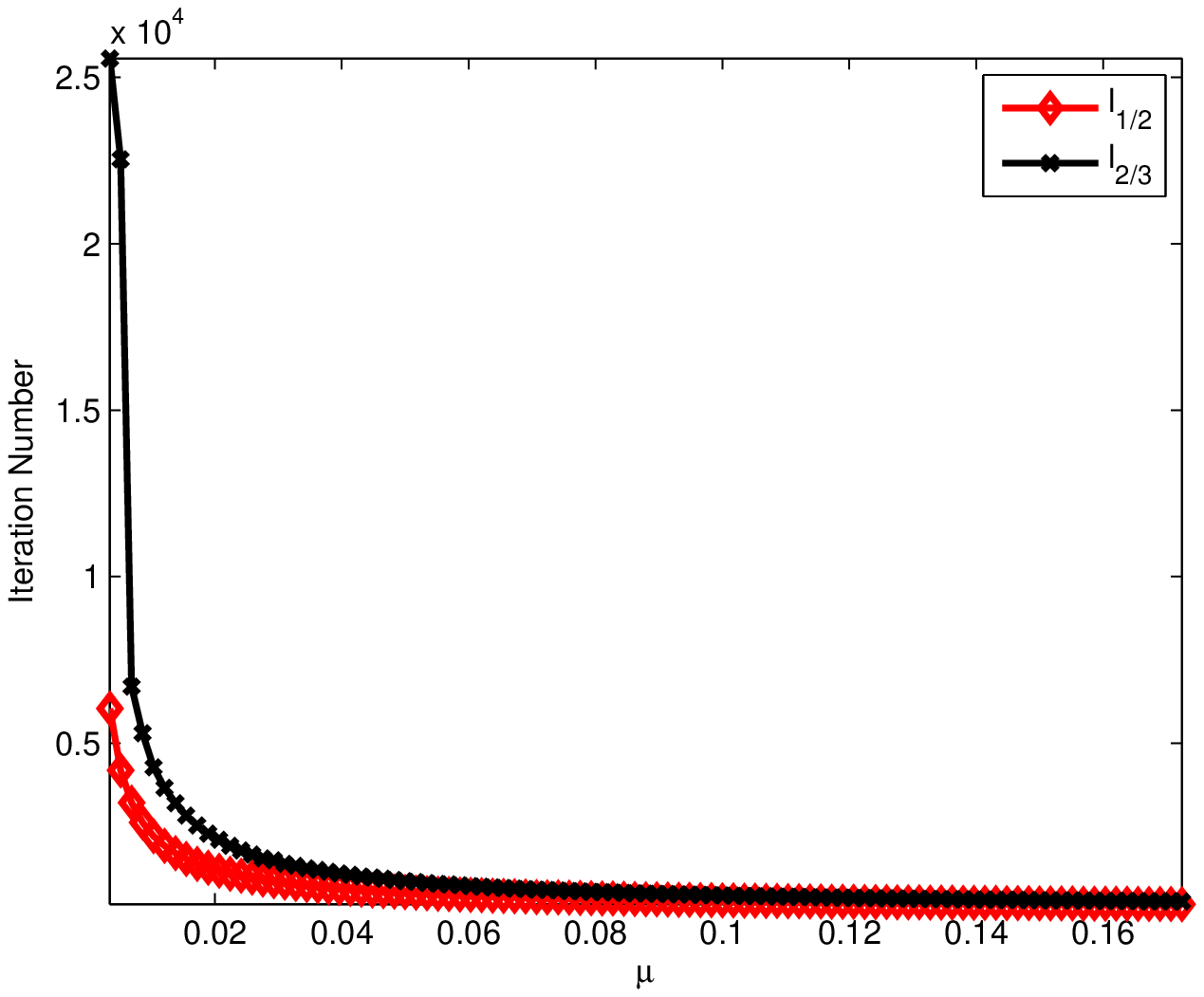}
\centerline{{\small (b) Iteration number}}
\end{minipage}
\begin{minipage}[b]{.32\linewidth}
\centering
\includegraphics*[scale=0.21]{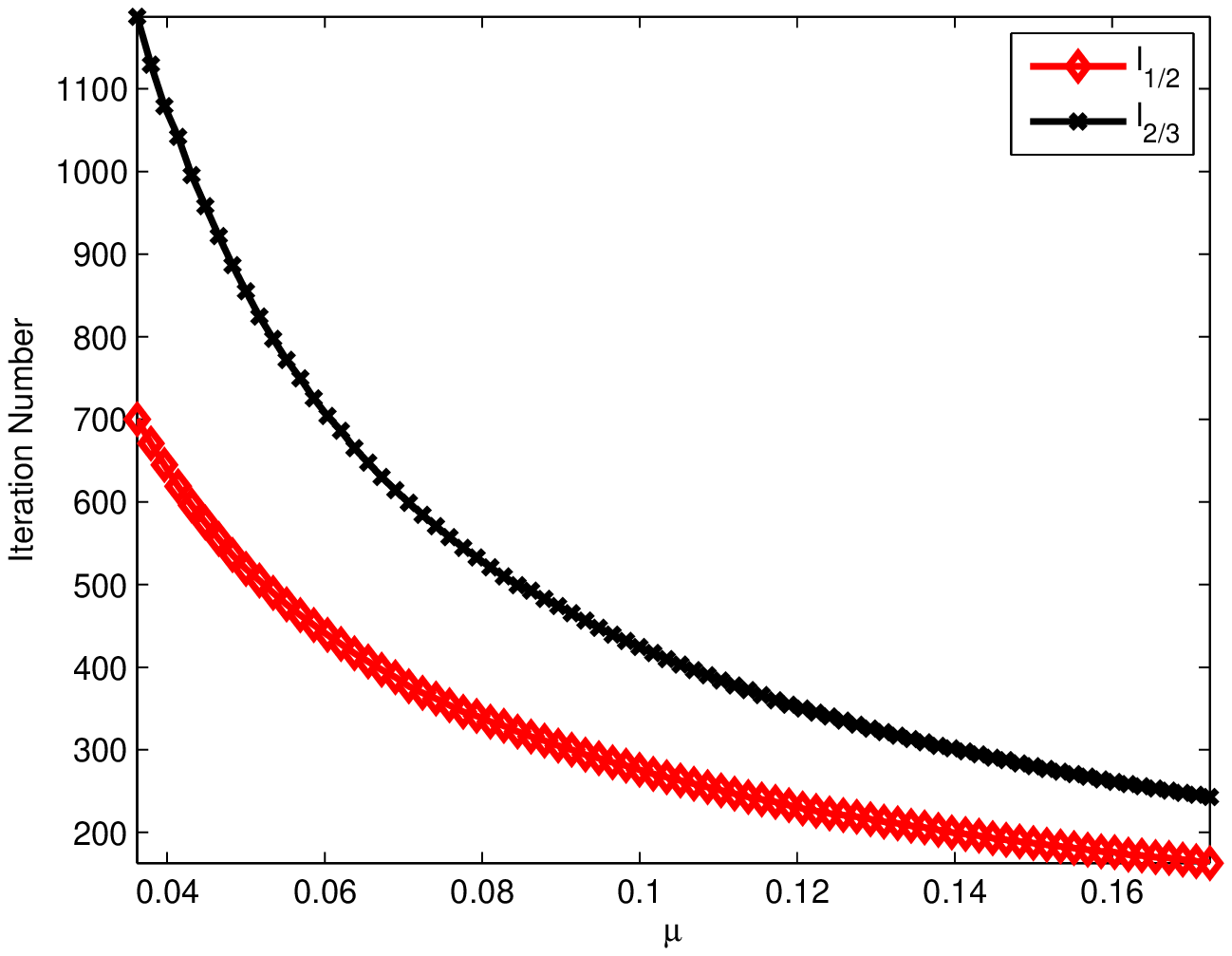}
\centerline{{\small (c) Detail}}
\end{minipage}
\caption{ Experiment for the effect of $\mu$.
(a) The trend of the recovery error.
(b) The trend of the required iteration numbers to achieve the setting accuracy.
(c) The detail trend of the required iteration numbers.
The regularization parameter $\lambda$ was taken as $0.001$,
the initialization was taken as the solution of the $l_1$-minimization problem and
the terminal rule of IJT algorithm was set as $\|x^{(n+1)}-x^{(n)}\|_2/\|x^{(n+1)}\|_2<10^{-10}$ for both penalties.
}
\label{Fig_mu}
\end{figure}

As shown by the iterative form (\ref{GGPA}) of IJT algorithm,  the
step size parameter $\mu$ is a crucial parameter of IJT algorithm.
In this subsection, we conducted a series of experiments to verify
the effect of $\mu$ on both the recovery precision and convergence
speed.
The measurement matrix and the true sparse signal were set the same
as in Subsection 6.1. We applied IJT algorithm for both $\phi(|z|) =
|z|^{1/2}$ and $\phi(|z|) = |z|^{2/3}$  with different $\mu$ to
recover the sparse signal from the given measurements. We varied
$\mu$ uniformly in the interval $(0,\|A\|_2^{-2})$ for 100 times.
The experimental results are shown in Fig. {\ref{Fig_mu}}.

From Fig. {\ref{Fig_mu}}(a), we can observe that $\mu$ has almost no effect on
the recovery quality of IJT algorithm for both penalties. While
the number of iterations required to attain the same terminal rule
decreases monotonically as $\mu$ increasing as demonstrated by Fig.
{\ref{Fig_mu}}(b) and (c). This phenomenon coincides with the common sense.
It demonstrates that when $\mu$ is larger, the algorithm converges
faster, and thus fewer iterations are required to attain a given
precision. More specifically, as shown by Fig. {\ref{Fig_mu}}(b), the number of
iterations decreases much sharper when $\mu<0.02$. Accordingly, we
recommend that in practical application of IJT algorithm, a larger
step size $\mu$ should be taken. In addition, we found that the
performance of IJT algorithm for $l_{1/2}$ regularization is
slightly better than the performance for $l_{2/3}$
regularization in the perspectives of both recovery quality and
iteration number, as shown in Fig. {\ref{Fig_mu}}. The additional advantage  of
IJT algorithm for $l_{1/2}$ regularization in the
perspective of cpu time was also demonstrated in the next subsection
over IJT algorithm for $l_{2/3}$ regularization.

\subsection{Comparisons with Reweighted Techniques}

\begin{figure}[!t]
\begin{minipage}[b]{.49\linewidth}
\centering
\includegraphics*[scale=0.33]{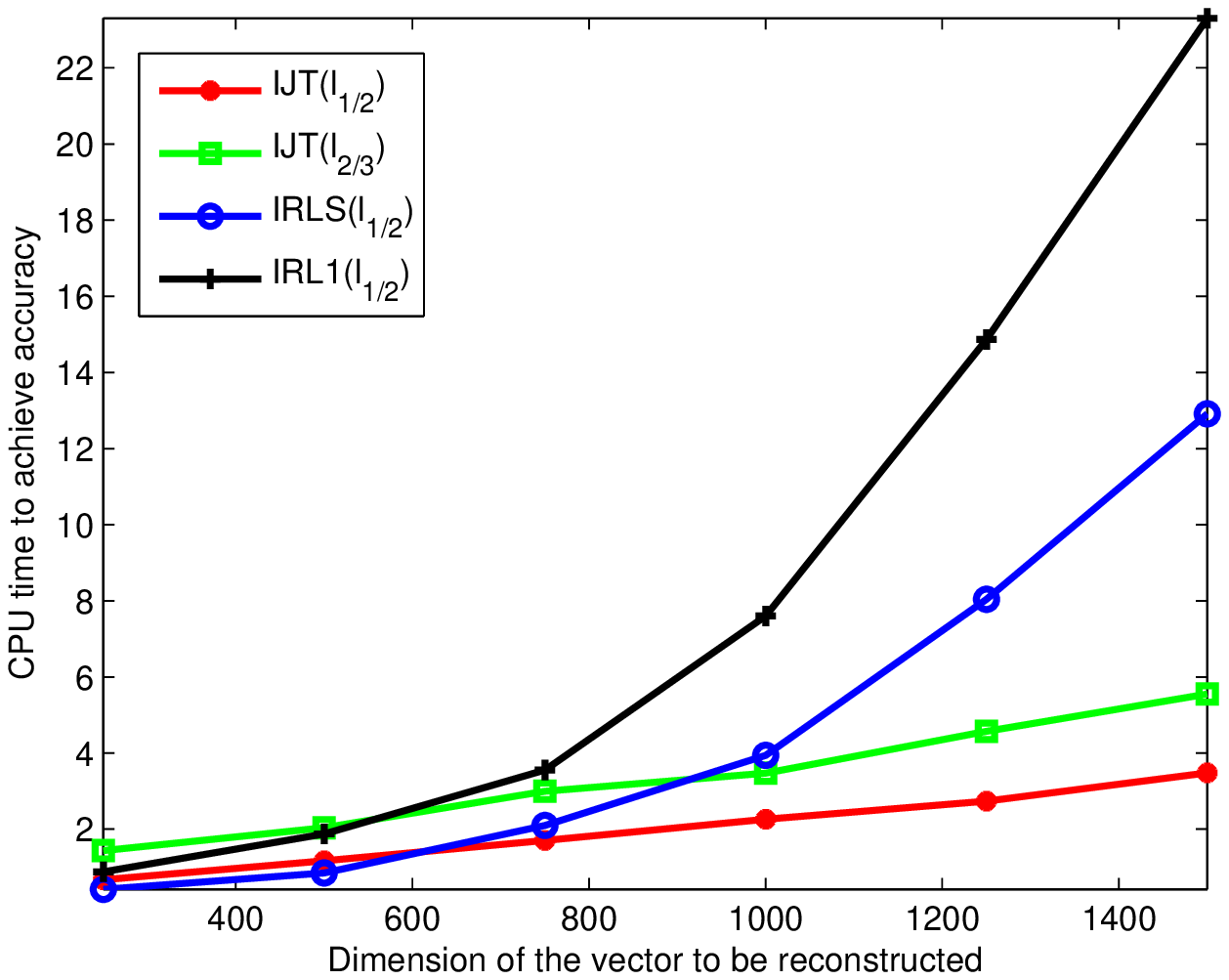}
\centerline{{\small (a) CPU time}}
\end{minipage}
\hfill
\begin{minipage}[b]{.49\linewidth}
\centering
\includegraphics*[scale=0.33]{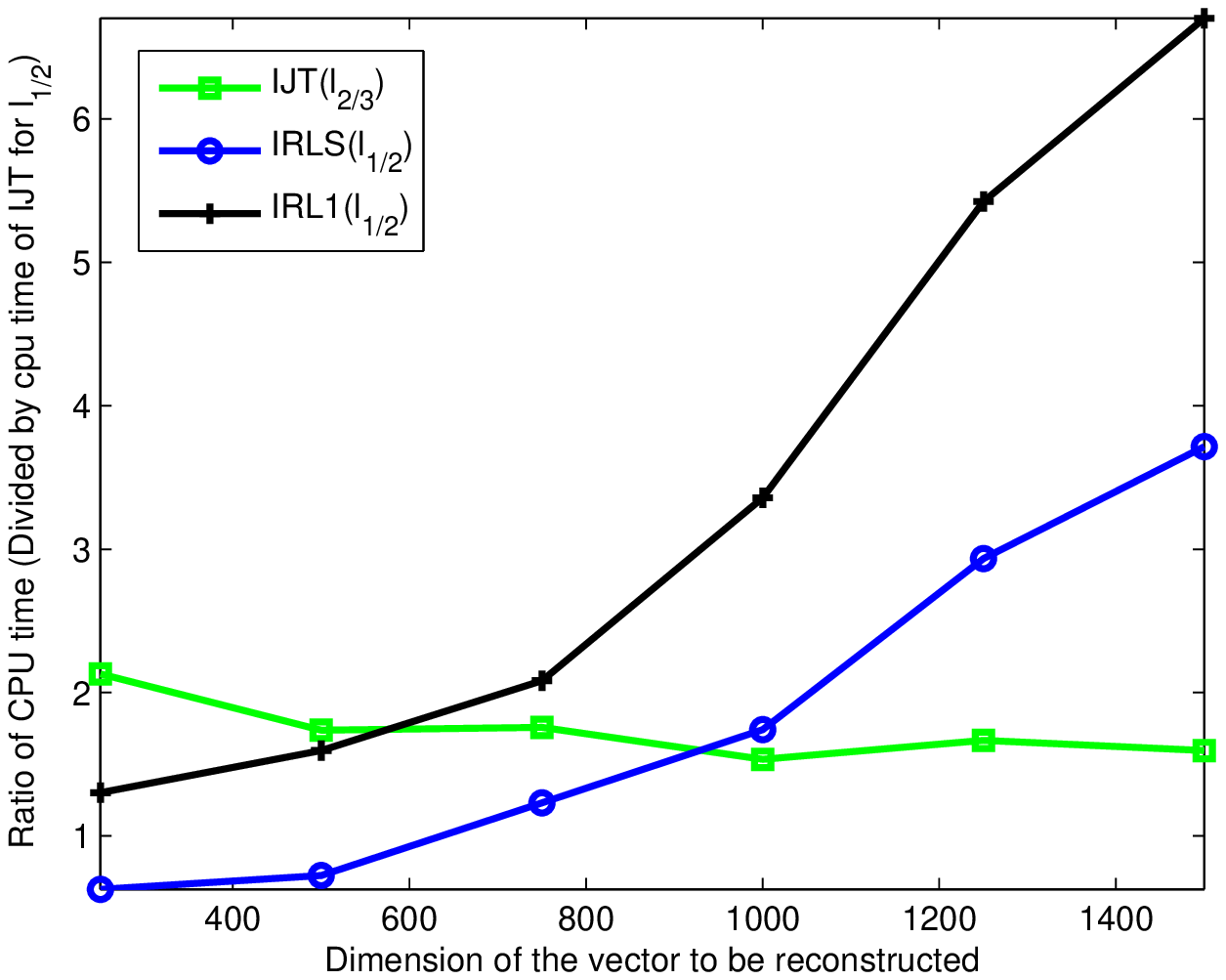}
\centerline{{\small (b) Ratio of CPU time}}
\end{minipage}
\hfill
\caption{ Experiment for comparison of CPU times of different algorithms including IJT, IRLS and IRL1 algorithms.
(a) The trends of CPU times of different algorithms.
(b) The trends of the ratios of CPU times (divided by the cpu time of IJT algorithm with $\phi(|z|) = |z|^{1/2}$).
}
\label{Fig_CompReweighted}
\end{figure}

This set of experiments were conducted to compare the time costs of
IJT algorithm, IRLS algorithm {\cite{DaubechiesIRLS}} and IRL1 algorithm {\cite{Candes2008RL1}} for solving
the same signal recovery problem with different settings $\{k,M,N\},$ where,
as in Subsection 8.2 in {\cite{DaubechiesIRLS}}, we took $k=5$, $N=\{250,500,750,1000,1250,1500\}$
and $M=N/5$.
We applied IJT algorithm for two different penalties, i.e., $\phi(|z|) = |z|^{1/2}$ and $\phi(|z|) = |z|^{2/3}$.
We implemented all algorithms using Matlab without any specific optimization.
In particular, we used the CVX Matlab package by
Michael Grant and Stephen Boyd (http://www.stanford.edu/ $\sim$boyd/cvx/) to
perform the weighted $l_{1}$-minimization at each iteration step of IRL1
algorithm. Again, the measurement matrix $A$ was taken to be the $M\times N$
dimensional matrices with i.i.d. Gaussian $\mathcal{N}(0,\frac{1}{M})$
entries. The experiment results are shown in Fig. {\ref{Fig_CompReweighted}}.
As shown in Fig. {\ref{Fig_CompReweighted}}(a), when $N$ is lower than $500$, IRLS
algorithm is slightly faster than IJT algorithm with $\phi(|z|)= |z|^{1/2}$.
This is due to that in the low-dimensional cases, the computational burden of
solving a low-dimensional least squares problem in IRLS is relatively low.
Nevertheless, when $N>500,$ it can be observed that IJT algorithm with $\phi(|z|)= |z|^{1/2}$
outperforms both IRLS and IRL1 algorithms in the perspective of CPU time. Furthermore,
we can observe from Fig. {\ref{Fig_CompReweighted}}(b) that as $N$ increases, the CPU times cost by IRL1
and IRLS algorithms increase much faster than IJT algorithm,
that is to say, the outperformance of IJT algorithm in time cost
can get more significant as dimension increases.


\section{Conclusion}

We have conducted a study of the convergence of IJT algorithm for a class of non-convex regularized optimization problems.
One of the most significant features of such class of iterative thresholding algorithms is that the associated thresholding functions are discontinuous with jump discontinuities.
Moreover, the corresponding thresholding functions are in general not nonexpansive due to the nonconvexity of the penalties.
Among such class of non-convex optimization problems, the $l_q$ ($0<q<1$) regularization problem is one of the most typical subclass.

The main contribution of this paper is the establishment of the convergence and rate-of-convergence results of IJT
algorithm for a certain class of non-convex optimization problems.
We first prove the finite support and sign convergence of IJT algorithm as long as $0<\mu<1/L,$ where $L$ is the Lipschitz constant of $\nabla F.$
Then we show the strong convergence of IJT algorithm under certain a rKL property.
Furthermore, we demonstrate that IJT algorithm converges to a local minimizer at an asymptotically linear rate under certain second-order conditions.
When applied to the $l_q$ regularization, IJT algorithm can converge to a local minimizer at an asymptotically linear rate
as long as the matrix satisfies a certain concentration property.
The obtained convergence
results to a local minimizer generalize those known for the \textit{soft} and \textit{hard} algorithms.
We have also provided a set of simulations to support the correctness of the
established theoretical assertions. The efficiency of IJT algorithm is
further compared through simulations with the known reweighted
techniques, another type of typical non-convex regularization algorithms.

\section*{Appendix}

\subsection{A non-KL function}
\begin{figure}[!t]
\begin{minipage}[b]{.49\linewidth}
\centering
\includegraphics*[scale=0.33]{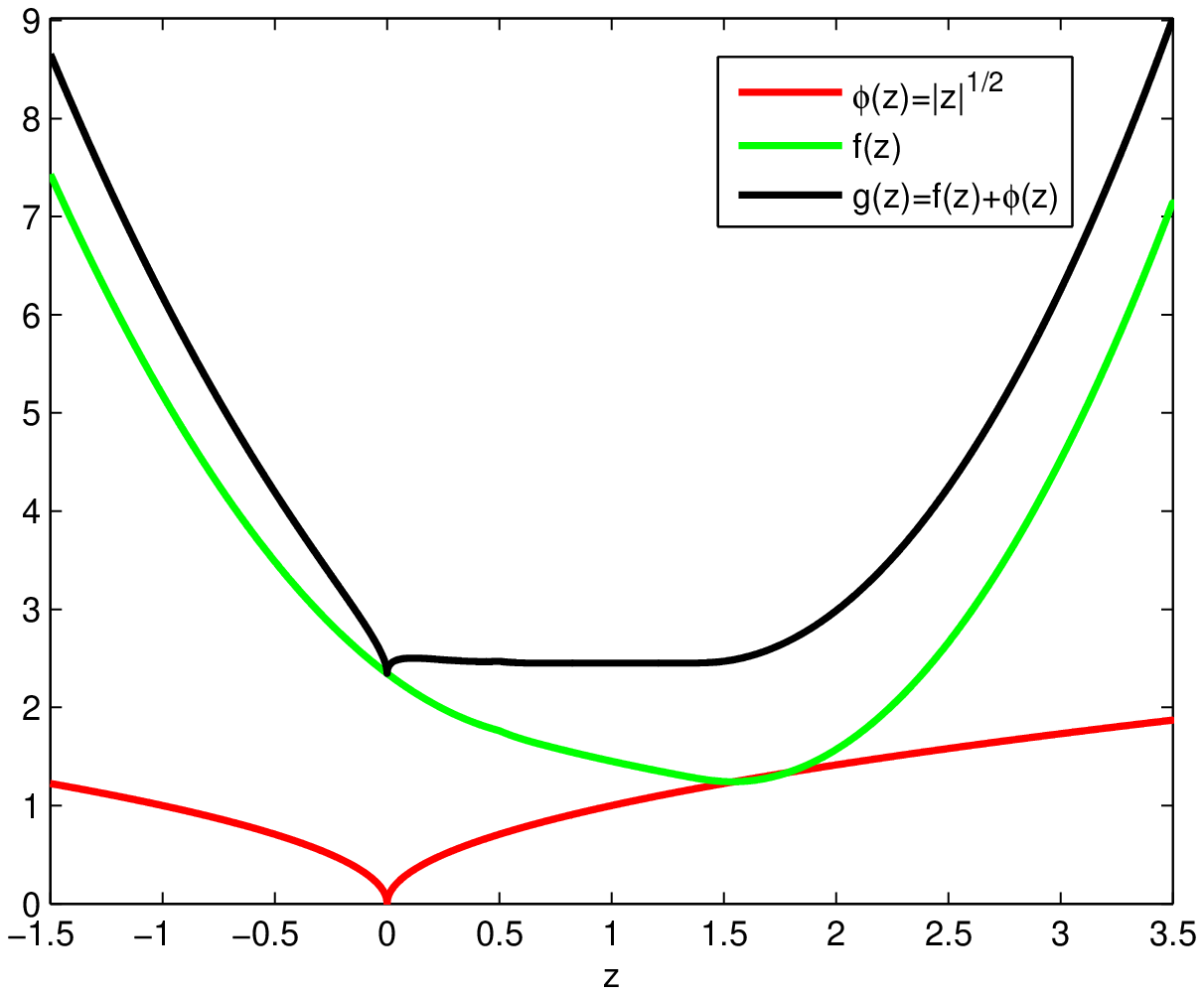}
\centerline{{\small (a) Figures of $\phi$, $f$ and $g$}}
\end{minipage}
\hfill
\begin{minipage}[b]{.49\linewidth}
\centering
\includegraphics*[scale=0.33]{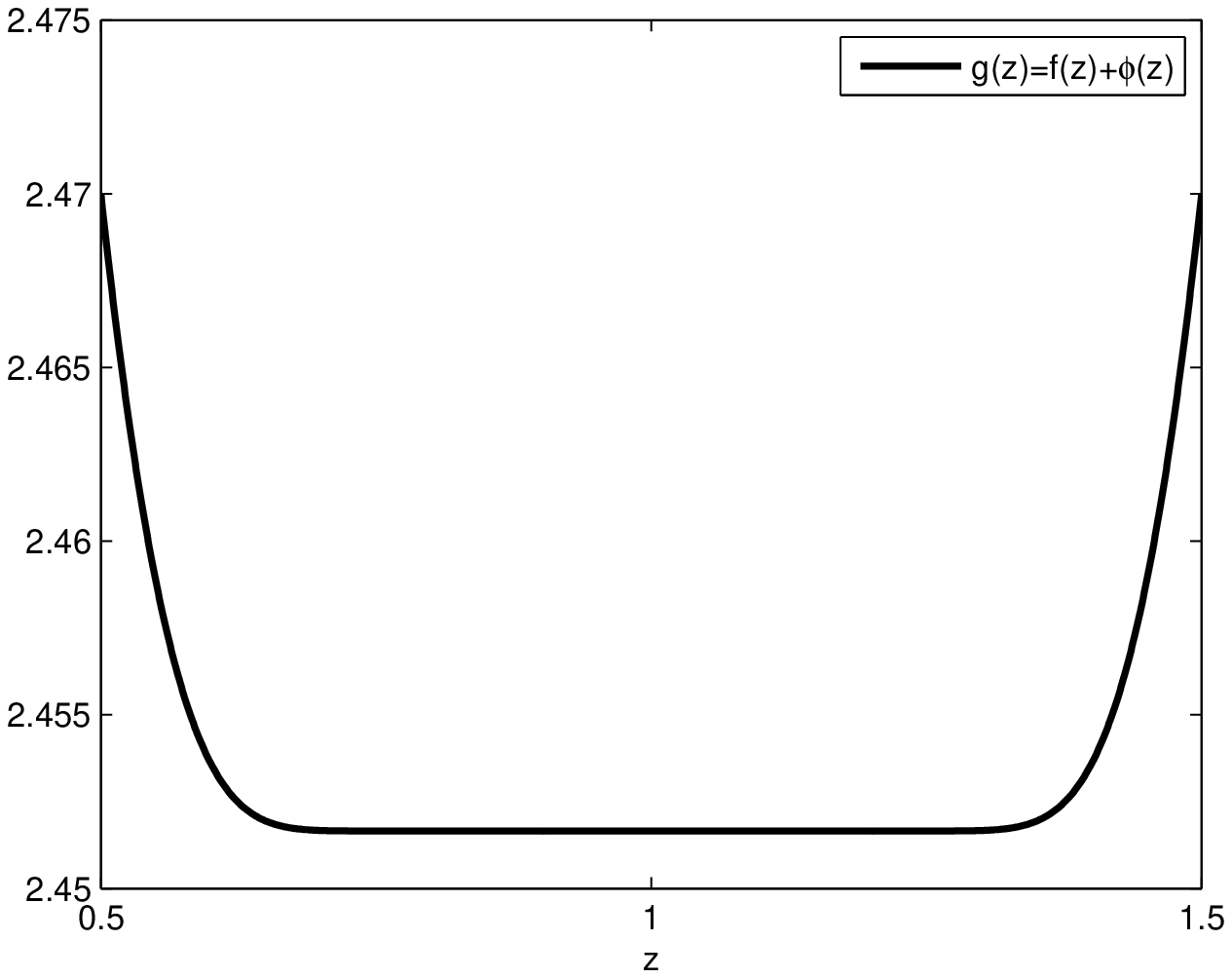}
\centerline{{\small (b) Detail figure of $g$}}
\end{minipage}
\hfill
\caption{ A specific function $g$ that is not KL function but satisfies Assumptions 1 and 2.
In this case, $\phi(|z|) = |z|^{1/2}$, $f$ is specified as in (\ref{SpecificCasef}) and $g=f+\phi$.
}
\label{Fig_nonKLFun}
\end{figure}

In the following, we give a specific one-dimensional function that satisfies Assumptions 1 and 2, but not a KL function.
Given any function $\phi$ satisfying Assumption 2, let $g=f+\phi$ with $f$ being defined as follows
\begin{equation}
f(z)=\left\{
\begin{array}
[c]{ll}%
a_1(z-b_1)^2+c_1, & \mbox{for} \ z\leq 1/2\\
\exp\left(-\frac{1}{(z-1)^2}\right)-\phi(z)+C, & \mbox{for} \ 1/2<z<1\\
C-\phi(1), & \mbox{for} \ z=1\\
\exp\left(-\frac{1}{(z-1)^2}\right)-\phi(z)+C, & \mbox{for} \ 1<z<3/2\\
a_2(z-b_2)^2+c_1, & \mbox{for} \  z\geq 3/2
\end{array}
\right.,
\label{SpecificCasef}
\end{equation}
where $e=exp(1)$,
$a_1 = 80e^{-4} - \frac{1}{2}\phi''(\frac{1}{2}),
b_1 = \frac{1}{2}+\frac{16e^{-4}+\phi'(\frac{1}{2})}{160e^{-4}-\phi''(\frac{1}{2})}
,$
$a_2 = 80e^{-4} - \frac{1}{2}\phi''(3/2),
b_2 = \frac{3}{2} - \frac{16e^{-4}-\phi'(\frac{3}{2})}{160e^{-4}-\phi''(\frac{3}{2})},$
$C = \phi(\frac{3}{2})+\max\left\{\phi(\frac{1}{2})+a_1(\frac{1}{2}-b_1)^2, \phi(\frac{3}{2})+a_2(\frac{3}{2}-b_2)^2 \right\}$,
$c_1 = C + e^{-4} -\phi(\frac{1}{2}) - a_1(\frac{1}{2}-b_1)^2,$
and
$ c_2 = C+e^{-4}-\phi(\frac{3}{2}) - a_2(\frac{3}{2}-b_2)^2.$
Thus,
\begin{equation}
g(z)=\left\{
\begin{array}
[c]{ll}%
a_1(z-b_1)^2+c_1 + \phi(|z|), & \mbox{for} \ z\leq 1/2\\
\exp\left(-\frac{1}{(z-1)^2}\right)+C, & \mbox{for} \ 1/2<z<1\\
C, & \mbox{for} \ z=1\\
\exp\left(-\frac{1}{(z-1)^2}\right)+C, & \mbox{for} \ 1<z<3/2\\
a_2(z-b_2)^2+c_1 + \phi(z), & \mbox{for} \  z\geq 3/2
\end{array}
\right..
\label{SpecificCaseg}
\end{equation}
When $1/2<z<3/2$, we define a function $h(z)$ as
\[
h(z)=\left\{
\begin{array}
[c]{ll}%
\exp\left(-\frac{1}{(z-1)^2}\right), & \mbox{for} \ 1/2<z<1\\
0, & \mbox{for} \ z=1\\
\exp\left(-\frac{1}{(z-1)^2}\right), & \mbox{for} \ 1<z<3/2\\
\end{array}
\right..
\label{SpecificCaseh}
\]
It can be easily checked that $f$ satisfies Assumption 1 due to the function $h$ is ${\cal{C}}^{\infty}$ and $\phi$ is ${\cal{C}}^{2}$ in the interval $(1/2,3/2)$.
However, according to {\cite{KLBolte2006}} (Sec. 1, page 1), it shows that $h$ fails to satisfy the KL inequality (\ref{KLIneq}) at $z=1$.
Therefore, $g$ must be not a KL function.
The figures of $f$ and $g$ are shown in Fig. {\ref{Fig_nonKLFun}} with $\phi(|z|) = |z|^{1/2}$.

\subsection{Proof of Lemma {\ref{Lemm_SuffCond_rKL}}}

\begin{proof}
Note that $z^*$ is a stationary point of $g$, i.e., $ \nabla g(z^*)=0$, then
\begin{align}
&|g(z) - g(z^*)|  =  |g(z) - g(z^*) - \nabla g(z^*)^T (z-z^*)| \nonumber\\
&\leq \int_{0}^1 \|\nabla g(z^*+t(z-z^*))-\nabla g(z^*)\|_2 \|z-z^*\|_2 dt.
\label{diffFun}
\end{align}
Since $g$ is twice continuously differentiable at $B(z^*,\epsilon_0)$,
then it obviously exists constants $L_g>0$ such that
$$
\|\nabla g(z^*+t(z-z^*))-\nabla g(z^*)\|_2 \leq L_g t \|z-z^*\|_2,
$$
for any $z\in B(z^*,\epsilon_0)$ and $t\in (0,1)$.
Thus, it follows
\begin{equation}
|g(z) - g(z^*)| \leq \frac{L_g}{2} \|z-z^*\|_2^2, \forall z\in B(z^*,\epsilon_0).
\label{diffFun1}
\end{equation}

On the other hand, for any $z\in B(z^*,\epsilon_0)$,
there exists a $t_0 \in (0,1)$ such that
\begin{align}
\label{DT1}
\|\nabla g(z)\|_2 = \|\nabla g(z) - \nabla g(z^*)\|_2 \\\nonumber
= \|\nabla^2 g(z^* + t_0(z-z^*)) (z-z^*)\|_2.
\end{align}
Since $\nabla^2 g(z^*)$ is nonsingular and by the continuity of $\nabla^2 g(z)$ at $B(z^*,\epsilon_0)$,
then there exists $0<\epsilon<\epsilon_0$ such that for any $z\in B(z^*,\epsilon),$
\[
\sigma_{\min}(\nabla^2 g(z^* + t_0(z-z^*))) \geq \min_{z\in B(z^*,\epsilon)} \sigma_{\min}(\nabla^2 g(z))>0.
\]
Denote $\sigma_{\epsilon,z^*} = \min_{z\in B(z^*,\epsilon)} \sigma_{\min}(\nabla^2 g(z)),$
then (\ref{DT1}) becomes
\begin{align}
\|\nabla g(z)\|_2 \geq \sigma_{\epsilon,z^*} \|z-z^*\|_2.
\label{DT2}
\end{align}
Let
$
C^* = \frac{L_g}{2\sigma_{\epsilon,z^*}^2}.
$
Combining (\ref{diffFun1}) and (\ref{DT2}), it implies
\[
|g(z) - g(z^*)| \leq C^* \|\nabla g(z)\|_2^2.
\]
Thus, we complete the proof of the lemma.

\end{proof}

\subsection{Proof of Lemma {\ref{Lemm_SuppConv}}}

\begin{proof}
\textbf{(i)}
By Property {\ref{Prop_SuffDecrease}}(b), there exists a sufficiently large positive integer $n_0$ such that $\|x^{n} - x^{n+1}\|_2 < \eta_{\mu}$ when $n>n_0$.
We first show that
\begin{equation}
I^{n+1} = I^n, \forall n>n_0
\label{Lemma3a}
\end{equation}
by contradiction.
Assume this is not the case, that is, $I^{n_1+1} \neq I^{n_1}$ for some $n_1>n_0$.
Then it is easy to derive a contradiction through
distinguishing the following two possible cases:

\textit{Case 1:} $I^{n_1+1}\neq I^{n_1}$ and $(I^{n_1+1}\cap I^{n_1}) \subset I^{n_1+1}.$
In this case, there exists an $i_{n_1}$ such that $i_{n_1}\in I^{n_1+1}\setminus I^{n_1}$.
By Lemma {\ref{Lemm_JumpThresFun}}, it then implies
\[
\Vert x^{n_1+1}-x^{n_1}\Vert_{2}\geq|x_{i_{n_1}}^{n_1+1}|\geq\min_{i\in I^{n_1+1}%
}|x_{i}^{n_1+1}|\geq \eta_{\mu},
\]
which contradicts to $\Vert x^{n_1+1}-x^{n_1}\Vert_{2}<\eta_{\mu}.$

\textit{Case 2:} $I^{n_1+1}\neq I^{n_1}$ and $(I^{n_1+1}\cap I^{n_1}) =I^{n_1+1}.$ Under
this circumstance, it is obvious that $I^{n_1+1}\subset I^{n_1}$.
Thus, there exists an $k_{n_1}$ such that $k_{n_1} \in I^{n_1} \setminus I^{n_1+1}$.
It then follows from Lemma {\ref{Lemm_JumpThresFun}} that
\[
\Vert x^{n_1+1}-x^{n_1}\Vert_{2}\geq|x_{k_{n_1}}^{n_1}|\geq\min_{i\in I^{n_1}}|x_{i}^{n_1}|\geq \eta_{\mu},
\]
and it contradicts to $\Vert x^{n_1+1}-x^{n_1}\Vert_{2}<\eta_{\mu}$.
Thus, (\ref{Lemma3a}) holds true.
It also means that the support set sequence $\{I^n\}$ converges.
We denote $I$ the limit of $I^n$.
Then for any $n>n_0,$ $I^n = I$.

\textbf{(ii)} For any limit point $x^* \in \cal{X}$,
there exits a subsequence $\{x^{n_j}\}$ converging to $x^*$, i.e.,
\begin{equation}
x^{n_j} \rightarrow x^*\ \ \text{as}\ \ j\rightarrow \infty.
\end{equation}
Thus, there exists a sufficiently large positive integer $j_0$ such that
$n_{j_0} > n_0$ and
$\|x^{n_j} - x^*\|_2 < \eta_{\mu}$ when $j\geq j_0$.
Similar to the proof procedure (i),
it can be also claimed that $I^{n_j}=Supp(x^*)$ for any $j\geq j_0$.
On the other hand, by (\ref{Lemma3a}), $I^{n_j} = I$.
Thus, for any limit point $x^*$, $Supp(x^*) = I$.

Taking $n^* = n_{j_0}$, then by the above analysis, it is obvious that the claims (a) and (b) in Lemma {\ref{Lemm_SuppConv}} hold true.

\textbf{(iii)}
As $I^{n}=I=Supp(x^*)$
for any $n>n^*$ and $x^*\in {\cal{X}}$,
it suffices to show that  $sign(x_i^{n+1}) = sign(x_i^n)$ and $sign(x_i^{n_j}) = sign(x_i^*)$ for any $i\in I$, $j\geq j_0$, $n>n^*$.
Similar to the first two parts of the proof, we will first check that $sign(x_i^{n+1}) = sign(x_i^n)$,
and then $sign(x_i^{n_j}) = sign(x_i^*)$ for any $i\in I$ by contradiction.
We now prove $sign(x_i^{n+1}) = sign(x_i^n)$ for any $i\in I$ and $n>n^*$.
Assume this is not the case. Then there exists an
$i^{\ast}\in I$ such that $sign(x_{i^*}^{n+1})\neq sign(x_{i^*}^{n})$, and
hence,
\[
sign(x_{i^{*}}^{n+1})sign(x_{i^{*}}^{n})=-1.
\]
From Lemma {\ref{Lemm_JumpThresFun}}, it is easy to check
\begin{align*}
\Vert x^{n+1}-x^{n}\Vert_{2}
& \geq|x_{i^{\ast}}^{n+1}-x_{i^{\ast}}^{n}|=|x_{i^{\ast}}^{n+1}|+|x_{i^{\ast}}^{n}|\\\nonumber
&\geq\min_{i\in I}\{|x_{i}^{n+1}|+|x_{i}^{n}|\}\geq 2\eta_{\mu},
\end{align*}
which contradicts again to $\Vert x^{n+1}-x^{n}\Vert_{2}<\eta_{\mu}$.
This contradiction shows
$sign(x^{n+1})=sign(x^{n})$ when $n>n^*$.
It follows that the sign sequence $\{sign(x^n)\}$ is convergent.
Let $S^*$ be the limit of the sign sequence $\{sign(x^n)\}$.
Similarly, we can also show that
$sign(x^{n_j}) = sign(x^*)$ whenever $j\geq j_0$.
Therefore, $sign(x^n)=S^*= sign(x^*)$ when $n>n^*$ and for any $x^*\in {\cal{X}}$.
This finishes the proof of Lemma {\ref{Lemm_SuppConv}}.
\end{proof}

\subsection{Proof of Theorem {\ref{Thm_ConvRate1}}}

\begin{proof}
Let
$C_1 = 1+\lambda \mu \phi''(e)$
and
$C_2 = \sqrt{1-2\mu \lambda_{\min}(\nabla_{II}^2 F(x^*)) + \mu^2L^2}.$
By the assumptions of Theorem {\ref{Thm_ConvRate1}}, it is easy to check that
$$C_1 > C_2>0.$$
Since both $c_{F}$ and $c_{\phi}$ approach to zero as $\varepsilon$ approaches zero,
then we can take a sufficiently small $0<\varepsilon<\eta_{\mu}$ such that
$$
0<c_F<\min\left\{\frac{(C_1-C_2)(C_1+3C_2)}{8\mu}, \lambda_{\min}(\nabla_{II}^2 F(x^*))\right\},
$$
and
$$
0<c_{\phi}<\frac{C_1-C_2}{2\lambda\mu}.
$$
Furthermore, let
$$
\alpha_{F,\varepsilon} = \lambda_{\min}(\nabla_{II}^2 F(x^*)) - c_F \ \text{and} \ \alpha_{\phi,\varepsilon} = -\phi''(e) + c_{\phi},
$$
then under assumptions of Theorem {\ref{Thm_ConvRate1}}, there hold
$0<\alpha_{F,\varepsilon}<L$ and $\alpha_{\phi,\varepsilon}>0,$
and further
\begin{align}
1-\lambda \mu \alpha_{\phi,\varepsilon}= 1+\lambda\mu \phi''(e) -\lambda\mu c_{\phi}  > \frac{C_1+C_2}{2}>0,
\label{rho*1}
\end{align}
\begin{align}
&1-2\mu \alpha_{F,\varepsilon} + \mu^2L^2  \geq 1-2\mu \alpha_{F,\varepsilon} + \mu^2\alpha_{F,\varepsilon}^2 \geq 0,
\label{rho*21}
\end{align}
\begin{align}
\label{rho*2}
&1-2\mu \alpha_{F,\varepsilon} + \mu^2L^2  = C_2^2 + 2\mu c_F \\\nonumber
&< C_2^2 + \frac{(C_1-C_2)(C_1+3C_2)}{4} = \left(\frac{C_1+C_2}{2}\right)^2.
\end{align}

Since $\{x^n\}$ converges to $x^*$,
then for any $0<\varepsilon < \eta_{\mu}$,
there exists a sufficiently large integer $n_0>n^*$ (where $n^*$ is specified as in Lemma {\ref{Lemm_SuppConv}}) such that
$$\|x^n-x^*\|_2< \varepsilon$$
when $n>n_0$. Let $I^n = Supp(x^n)$.
By Lemma {\ref{Lemm_SuppConv}}, it holds
$I^n = I$ and $sign(x^n) = sign(x^*)$ when $n>n_0$.
Furthermore, by Property {\ref{Prop_OptCond}}, for any $i\in I$,
$$
x_i^* + \lambda \mu sign(|x_i^*|)\phi'(|x_i^*|) = x_i^* - \mu [\nabla F(x^*)]_i,
$$
and
$$
x_i^{n+1} + \lambda \mu sign(|x_i^{n+1}|)\phi'(|x_i^{n+1}|) = x_i^n - \mu [\nabla F(x^n)]_i,
$$
when $n>n_0$.
Consequently,
\begin{align*}
&(x_I^{n+1} - x_I^*) + \lambda \mu (\phi_1(x_I^{n+1})-\phi_1(x_I^*))\\\nonumber
&= (x_I^{n} - x_I^*) - \mu ([\nabla F(x^n)]_I - [\nabla F(x^*)]_I),
\end{align*}
and then
\begin{align}
&\|x_I^{n+1} - x_I^*\|_2^2 + \lambda \mu \langle \phi_1(x_I^{n+1})-\phi_1(x_I^*), x_I^{n+1} - x_I^* \rangle = \nonumber\\
& \langle  x_I^{n+1}-x_I^*, (x_I^{n} - x_I^*) - \mu ([\nabla F(x^n)]_I - [\nabla F(x^*)]_I) \rangle .
\label{lemma4.1}
\end{align}
By (\ref{Condphi}), the left side of (\ref{lemma4.1}) satisfies
\begin{align*}
&\|x_I^{n+1} - x_I^*\|_2^2 + \lambda \mu \langle \phi_1(x_I^{n+1})-\phi_1(x_I^*), x_I^{n+1} - x_I^* \rangle \\\nonumber
&\geq (1-\lambda\mu \alpha_{\phi,\varepsilon}) \|x_I^{n+1} - x_I^*\|_2^2,
\end{align*}
and the right side of (\ref{lemma4.1}) satisfies
\begin{align*}
&\langle  x_I^{n+1}-x_I^*, (x_I^{n} - x_I^*) - \mu ([\nabla F(x^n)]_I - [\nabla F(x^*)]_I) \rangle \leq\nonumber\\
& \|x_I^{n+1}-x_I^*\|_2  \|(x_I^{n} - x_I^*) - \mu ([\nabla F(x^n)]_I - [\nabla F(x^*)]_I)\|_2.
\end{align*}
Without loss of generality, we assume that $\|x_I^{n+1} -x_I^*\|_2 >0$, otherwise, it demonstrates that IJT algorithm converges to $x^*$ in finite iterations.
Thus, it becomes
\begin{align}
\label{lemma4.2}
&(1-\lambda\mu \alpha_{\phi,\varepsilon}) \|x_I^{n+1} - x_I^*\|_2 \\\nonumber
&\leq \|(x_I^{n} - x_I^*) - \mu ([\nabla F(x^n)]_I - [\nabla F(x^*)]_I)\|_2.
\end{align}
Furthermore,  by (\ref{CondF}), it follows
\begin{align}
&\|(x_I^{n} - x_I^*) - \mu ([\nabla F(x^n)]_I - [\nabla F(x^*)]_I)\|_2^2 \nonumber\\
& = \|x_I^{n} - x_I^*\|_2^2 + \mu^2 \|[\nabla F(x^n)]_I - [\nabla F(x^*)]_I\|_2^2 \nonumber\\
& - 2\mu \langle x_I^{n} - x_I^*, [\nabla F(x^n)]_I - [\nabla F(x^*)]_I \rangle \nonumber\\
& \leq (1-2\mu \alpha_{F,\varepsilon} + \mu^2 L^2) \|x_I^n - x_I^*\|_2^2.
\label{lemma4.3}
\end{align}
Combing (\ref{lemma4.2}) and (\ref{lemma4.3}), it implies
\begin{align*}
\|x_I^{n+1} - x_I^*\|_2 \leq
\frac{\sqrt{1-2\mu \alpha_{F,\varepsilon} + \mu^2 L^2}}{1-\lambda\mu \alpha_{\phi,\varepsilon}} \|x_I^n - x_I^*\|_2.
\end{align*}
Let
$$
\rho^* = \frac{\sqrt{1-2\mu \alpha_{F,\varepsilon} + \mu^2 L^2}}{1-\lambda\mu \alpha_{\phi,\varepsilon}}.
$$
By (\ref{rho*1})-(\ref{rho*2}),
it is easy to check that
\begin{equation*}
0<\rho^*<1.
\end{equation*}
Thus, when $n>n_0$
\begin{align}
\label{LinConLemm4}
&\|x^{n+1} - x^*\|_2 = \|x_I^{n+1} - x_I^*\|_2\\\nonumber
& \leq \rho^* \|x_I^{n} - x_I^*\|_2 = \rho^* \|x^{n} - x^*\|_2.
\end{align}
Consequently, the asymptotic convergence rate of IJT algorithm is linear.

Moreover, the posteriori error bound can be easily derived by the triangle inequality
$$
\|x^{n} - x^*\|_2 \leq \|x^{n+1} - x^*\|_2 + \|x^{n+1} - x^n\|_2
$$
and (\ref{LinConLemm4}).
Therefore, we have completed the proof of Theorem {\ref{Thm_ConvRate1}}.
\end{proof}

\subsection{Proof of Theorem {\ref{Thm_ConvRate2}}}

\begin{proof}
Let
\begin{equation}
c_1 = \frac{1-\mu\lambda_{min}(\nabla_{II}^2F(x^*))}{1+\lambda\mu\phi''(e)}.
\label{Const_c1}
\end{equation}
By the assumptions of Theorem {\ref{Thm_ConvRate2}}, it holds $0<c_1<1$.
For any $0<c<1$, let
\begin{equation}
g(c) = \max_{i\in I} \max_{\{x_i: |x_i-x_i^*|<c\eta_{\mu}\}}
\left\{\frac{\lambda\mu|\phi'''(|x_i|)|}{2|1 + \lambda\mu \phi''(|x_i^*|)|}\right\},
\label{const_g}
\end{equation}
and
\begin{equation}
c_{\epsilon}(c) = \frac{1-c_1-\epsilon}{g(c)\eta_{\mu}},
\label{c_epsilon}
\end{equation}
for some $0<\epsilon < 1-c_1.$
Since $g(c)$ is non-decreasing with respective to $c$, and thus $c_{\epsilon}(c)$ is non-increasing with respect to $c$.
Therefore, there exists a positive constant $c^*$ such that
\begin{equation}
0<c^*<1 \ \text{and}\ c^*< c_{\epsilon}(c^*).
\label{Const_c*}
\end{equation}

Since $\{x^n\}$ converges to $x^*$, then there exists an $n^{**}>n^*$ (where $n^*$ is specified as in Lemma {\ref{Lemm_SuppConv}}), when $n>n^{**}$, it holds
$$
\|x^n-x^*\|_2<c^*\eta_{\mu}.
$$
By Lemma {\ref{Lemm_SuppConv}}, when $n>n^{**}$, it holds
$I^n = I$ and $sign(x^n) = sign(x^*)$ ,
and thus $\|x^n-x^*\|_2 = \|x_{I}^n - x_{I}^*\|_2$.
By Property {\ref{Prop_OptCond}}, for any $i\in I$,
\begin{align*}
&(x_i^n - x_i^*) - \mu ([\nabla F(x^n)]_i - [\nabla F(x^*)]_i) \nonumber\\
&= (x_{i}^{n+1} - x_{i}^{*}) + sign(x_{i}^{*}) \lambda \mu( \phi'(|x_{i}^{n+1}|)- \phi'(|x_{i}^{*}|)).
\end{align*}
By Taylor expansion, for any $i\in I$,
there exists an $\xi_i \in (0,1)$,
such that
\begin{align*}
& \phi'(|x_i^{n+1}|) - \phi'(|x_i^*|)= \nonumber\\
& sign(x_i^*) \phi''(|x_i^*|) (x_i^{n+1} - x_i^{*}) + \frac{1}{2} \phi'''(|x_i^{\xi}|)(x_i^{n+1} - x_i^{*})^2,
\label{Taylorphi}
\end{align*}
where $x_i^{\xi} = x_i^* + \xi_i (x_i^{n+1} - x_i^*)$.
Let $h^n = x^n - x^*$, then by the above two inequalities, it follows
\begin{equation}
\Lambda_1 h_{I}^{n+1} + \Lambda_2 (h_{I}^{n+1}\odot h_{I}^{n+1})=
 h_{I}^n - \mu ([\nabla F(x^n)]_{I} - [\nabla F(x^*)]_{I}),
\label{Th6Eq2}
\end{equation}
where $\odot$ denotes the Hadamard product or elementwise product,
$\Lambda_1$ and $\Lambda_2$ are two different diagonal matrices with
\begin{align}
\label{Lambda}
&\Lambda_1(i,i) = 1 + \lambda\mu \phi''(|x_i^*|), \\\nonumber
&\Lambda_2(i,i) = \frac{1}{2}sign(x_i^*)\lambda\mu \phi'''(x_i^{\xi}).
\end{align}

Moreover, by the twice differentiability of $F$ at $x^*$, we have
\begin{equation}
[\nabla F(x^n)]_{I} - [\nabla F(x^*)]_{I} = \nabla_{II}^2 F(x^*)h_{I}^n + o(\|h_{I}^n\|_2).
\label{TaylorF2}
\end{equation}
Plugging (\ref{TaylorF2}) into (\ref{Th6Eq2}), it becomes
\begin{equation*}
\Lambda_1 h_{I}^{n+1} + \Lambda_2 (h_{I}^{n+1}\odot h_{I}^{n+1}) = (\mathbf{I} - \mu \nabla_{II}^2 F(x^*)) h_{I}^n + o(\|h_{I}^n\|_2),
\label{Th6Eq3}
\end{equation*}
where $\mathbf{I}$ denotes as the identity matrix with the size $|I|\times |I|$ with $|I|$ being the cardinality of the set $I$.
By the assumptions of Theorem {\ref{Thm_ConvRate2}}, for any $i\in I$,
\begin{align*}
&\Lambda_1(i,i) = 1 + \lambda\mu \phi''(|x_i^*|) \\\nonumber
&\geq 1+\lambda\mu \phi''(e) > 1-\mu \lambda_{\min}(\nabla^2_{II} F(x^*)) \geq 0,
\end{align*}
thus, $\Lambda_1$ is invertible.
Then it follows
\begin{align}
\label{Th6Eq4}
& h_{I}^{n+1} = \Lambda_1^{-1} (\mathbf{I} - \mu \nabla_{II}^2 F(x^*)) h_{I}^n \\\nonumber
& -\Lambda_1^{-1} \Lambda_2 (h_{I}^{n+1}\odot h_{I}^{n+1}) + o(\|h_{I}^n\|_2).
\end{align}
By the definition of $o(\|h_{I}^n\|_2)$, there exists a constant $c_{\epsilon}^*$ (depending on $\epsilon$) such that
$$
|o(\|h_{I}^n\|_2)| \leq \epsilon \|h_{I}^n\|_2
$$
when $\|h_{I}^n\|_2<c_{\epsilon}^* \eta_{\mu}$.
Thus, we can take $c_0 = \min\{c^*, c_{\epsilon}^*\}<1$ and $n_0>n^{**}$ such that when $n>n_0$,
$$
\|x^n-x^*\|_2 <c_0\eta_{\mu}.
$$
Then (\ref{Th6Eq4}) implies that
\begin{align*}
\|h_{I}^{n+1} \|_2
&\leq \|\Lambda_1^{-1} (I - \mu \nabla_{II}^2 F(x^*)) h_{I}^n\|_2 \nonumber\\
&+ \epsilon \|h_{I}^n\|_2+ \|\Lambda_1^{-1} \Lambda_2 (h_{I}^{n+1}\odot h_{I}^{n+1})\|_2 \nonumber\\
&\leq \|\Lambda_1^{-1} (I - \mu \nabla_{II}^2 F(x^*))\|_2 \| h_{I}^n\|_2 \nonumber\\
&+ \epsilon \|h_{I}^n\|_2+ g(c^*) \| h_{I}^{n+1}\|_2^2 \nonumber \\
&\leq \left(\frac{1-\mu \lambda_{min}(\nabla_{II}^2 F(x^*))}{1+\lambda\mu \phi''(e)} + \epsilon\right)
\| h_{I}^n\|_2 \nonumber\\
&+ g(c^*) \| h_{I}^{n+1}\|_2^2\nonumber\\
& \leq (c_1 + \epsilon)\|h_{I}^n\|_2 + g(c^*) c^* \eta_{\mu}\| h_{I}^{n+1}\|_2,
\label{Th6Eq5}
\end{align*}
where the second inequality holds for the definition of $g(c^*)$ as specified in (\ref{const_g}) and $c^* \geq c_0$,
the third inequality holds for $\lambda_{\max}(I - \mu \nabla_{II}^2 F(x^*)) \leq 1-\mu \lambda_{min}(\nabla_{II}^2 F(x^*))$ and $\min_{i\in I} |\Lambda_1(i,i)| \geq 1+\lambda\mu \phi''(e)>0$, the last inequality holds for $\|h_I^{n+1}\|_2<c^*\eta_{\mu}$ and the definition of $c_1$ as specified in (\ref{Const_c1}).
Furthermore, by (\ref{c_epsilon}) and (\ref{Const_c*}), it holds
$$
1-c^*g(c^*)\eta_{\mu} >c_1+\varepsilon>0.
$$
Therefore, it implies that
\begin{equation*}
\|h_{I}^{n+1}\|_2 \leq \frac{c_1+\epsilon}{1-c^*g(c^*)\eta_{\mu}}{\|h_{I}^n\|_2},
\label{LinCon1}
\end{equation*}
and then
\begin{equation*}
\|x^{n+1}-x^*\|_2 \leq  \frac{c_1+\epsilon}{1-c^*g(c^*)\eta_{\mu}}{\|x^n-x^*\|_2}.
\label{LinCon2}
\end{equation*}
Let
$
\rho = \frac{c_1+\epsilon}{1-c^*g(c^*)\eta_{\mu}},
$
then
$
0<\rho<1.
$
Thus, the asymptotic convergence rate of IJT algorithm is linear.

Moreover, the error bound can be easily derived by the asymptotic convergence rate and the triangle inequality.

\end{proof}

\end{document}